\font\fiverm=cmr5

\font\sixrm=cmr6

\font\sevenrm=cmr7

\font\eightrm=cmr8

\font\ninerm=cmr9

\font\bigrm=cmr10 scaled \magstep1

\font\eightbf=cmbx8
\font\ninebf=cmbx9
\font\bigbf=cmbx10 scaled \magstep1
\font\Bigbf=cmbx10 scaled \magstep2

\font\tengoth=eufm10
\font\sevengoth=eufm7
\font\fivegoth=eufm5

 \newfam\gothfam
\textfont \gothfam=\tengoth
\scriptfont \gothfam=\sevengoth \scriptscriptfont \gothfam=\fivegoth


\newfam\srmfam \textfont
\srmfam=\eightrm \scriptfont
 \srmfam=\sixrm \scriptscriptfont \srmfam=\fiverm

\font\tengothb=eufb10
\font\sevengothb=eufb7
\font\fivegothb=eufb5 \newfam\gothbfam \textfont \gothbfam=\tengothb
\scriptfont
\gothbfam=\sevengothb \scriptscriptfont \gothbfam=\fivegothb

\font\tenmath=msbm10
\font\sevenmath=msbm7
\font\fivemath=msbm5 \newfam\mathfam \textfont \mathfam=
\tenmath \scriptfont
\mathfam=\sevenmath \scriptscriptfont \mathfam=\fivemath
\def\math{\fam\mathfam\tenmath}

\parindent=0mm

\def\titre#1{
\centerline{\Bigbf #1}\nobreak\nobreak\vglue 10mm\nobreak}

\def
\paragraphe#1{\bigskip\goodbreak {\bigbf #1}\nobreak\vglue 12pt\nobreak}
\def\alinea#1{\medskip\allowbreak{\bf#1}\nobreak\vglue 9pt\nobreak}
\def\ssq{\smallskip\qquad}
\def\msq{\medskip\qquad}


\def\th#1{\bigskip\goodbreak {\bf Th\'eor\`eme #1.} \par\nobreak \sl }
\def\prop#1{\bigskip\goodbreak {\bf Proposition #1.} \par\nobreak \sl }
\def\lemme#1{\bigskip\goodbreak {\bf Lemme #1.} \par\nobreak \sl }
\def\cor#1{\bigskip\goodbreak {\bf Corollaire #1.} \par\nobreak \sl }
\def\dem{\bigskip\goodbreak \it D\'emonstration. \rm}
\def\ndem{\bigskip\goodbreak \rm}
\def\qed{\par\nobreak\hfill $\bullet$ \par\goodbreak}


\let\wt=\widetilde
\def\uple#1#2{#1_1,\ldots ,{#1}_{#2}}
\def\corde#1#2-#3{{#1}_{#2},\ldots ,{#1}_{#3}}

\def\ordcorde#1#2-#3{{#1}_{#2} \le \cdots \le {#1}_{#3}}
\def\strictordcorde#1#2-#3{{#1}_{#2} < \cdots < {#1}_{#3}}
\def \restr#1{\mathstrut_{\textstyle |}\raise-6pt\hbox{$\scriptstyle #1$} }
\def \srestr#1{\mathstrut_{\scriptstyle |}\hbox to -1.5pt{}
\raise-4pt\hbox{$\scriptscriptstyle #1$}}
\def \inver{^{-1}}
\def\dbar{d\!\!\hbox to 4.5pt{\hfill\vrule height 5.5pt
depth -5.3pt width 3.5pt}}

\def\frac#1#2{{\textstyle {#1\over #2}}}

\def\R{{\math R}}
\def\C{{\math C}}

\def\N{{\math N}}

\def\I{{\math I}}

\def\fleche#1{\mathop{\hbox to #1 mm{\rightarrowfill}}\limits}
\def\gfleche#1{\mathop{\hbox to #1 mm{\leftarrowfill}}\limits}
\def\inj#1{\mathop{\hbox to #1 mm{$\lhook\joinrel$\rightarrowfill}}
\limits}
\def\ginj#1{\mathop{\hbox to #1 mm{\leftarrowfill$\joinrel\rhook$}}\limits}
\def\surj#1{\mathop{\hbox to #1 mm{\rightarrowfill\hskip 2pt\llap
{$\right arrow$}}}\limits}
\def\gsurj#1{\mathop{\hbox to #1 mm{\rlap{$\leftarrow$}
\hskip 2pt \leftar rowfill}}\limits}

\let\O = \Omega

\let \t = \theta
\let \l = \lambda

\let \f =\varphi

\let\p  =\partial


\def \g#1{\hbox{\tengoth #1}}
\def \sg#1{\hbox{\sevengoth #1}}

\def\Cal #1{{\cal #1}}


\def \ad{\mathop{\hbox{\rm ad}}\nolimits}

\def \mop#1{\mathop{\hbox{\rm #1}}\nolimits}
\def \smop#1{\mathop{\hbox{\sevenrm #1}}\nolimits}

\def \mopl#1{\mathop{\hbox{\rm #1}}\limits}


\def\bib #1{\null\medskip \strut\llap{[#1]\quad}}
\def\cite#1{[#1]} \magnification=\magstep1 \parindent=0cm
\def\titre#1{
\centerline{\Bigbf #1}
\vskip 16pt}
\def
\paragraphe#1{\bigskip {\bigbf #1}
\vskip 12pt}
\def\alinea#1{\medskip{\bf #1}
\vskip 6pt}
\def\ssq{\smallskip\qquad}
\def\msq{\medskip\qquad}

\def\I{\hbox{\rm I}}
\def\II{\hbox{\rm II}}
\def\III{\hbox{\rm III}}
\def\IV{\hbox{\rm IV}}
 \null
\long\def\ignore#1{}
\overfullrule=0cm
\centerline{\bigrm Ali Baklouti\footnote{\tenrm *}
{\eightrm Facult\'e des Sciences de Sfax,
D\'epartement de Math\'ematiques, 3038 Sfax, Tunisie.
Ali.Baklouti@fss.rnu.tn},
\bigrm Sami Dhieb\footnote{\tenrm **}
{\eightrm Facult\'e des Sciences de Sfax,
D\'epartement de Math\'ematiques, 3038 Sfax, Tunisie.
Sami.Dhieb@fss.rnu.tn},
\bigrm Dominique Manchon\footnote{\tenrm ***}{\eightrm CNRS - UMR 6620,
Clermont-Ferrand, France. manchon@math.univ-bpclermont.fr}}
\medskip \titre{Orbites coadjointes et vari\'et\'es caract\'eristiques}
\hfill {\sl A la m\'emoire de N.V. Pedersen \footnote{}{}\footnote{}{\eightbf Recherche
support\'ee par l'action d'\'echange CNRS-DGRSRT 01/R 15 04.}}
\vskip 15mm { \baselineskip=10pt \ninebf Abstract \ninerm : The
purpose of the present work is to describe a dequantization
procedure for topological modules over a deformed algebra. We
define the characteristic variety of a topological module as the
common zeroes of the annihilator of the representation obtained by
setting the deformation parameter to zero. On the other hand, the
Poisson characteristic variety is defined as the common zeroes of
the ideal obtained by considering the annihilator of the deformed
representation, and only then setting the deformation parameter to
zero. \vskip 2mm\qquad Using Gabber's theorem, we show the
involutivity of the characteristic variety. The Poisson
characteristic variety is indeed a Poisson subvariety of the
underlying Poisson manifold. We compute explicitly the
characteristic variety in several examples in the Poisson-linear
case, including the dual of any exponential solvable Lie algebra.
In the nilpotent case, we show that any coadjoint orbit appears as
the Poisson characteristic variety of a well chosen topological
module.
\vskip 12mm
\ninebf R\'esum\'e \ninerm : Nous pr\'esentons
dans ce travail un proc\'ed\'e de d\'equantification pour des modules
topologiques sur une alg\`ebre d\'eform\'ee. Nous d\'efinissons la
vari\'et\'e caract\'eristique d'un module topologique comme
l'ensemble des z\'eros communs de l'annulateur de la
repr\'esentation obtenue en annulant le param\`etre de
d\'eformation. Nous d\'efinissons par ailleurs la vari\'et\'e de
Poisson caract\'eristique comme l'ensemble des z\'eros communs de
l'id\'eal obtenu par quotient en annulant le param\`etre de
d\'eformation dans l'annulateur de la repr\'esentation
d\'eform\'ee. \vskip 2mm\qquad Nous montrons \`a l'aide du
th\'eor\`eme de Gabber l'involutivit\'e de la vari\'et\'e
caract\'eristique. La vari\'et\'e de Poisson caract\'eristique est
une sous-vari\'et\'e de Poisson de la vari\'et\'e de Poisson
sous-jacente. Nous explicitons la vari\'et\'e caract\'eristique
dans plusieurs exemples, incluant le dual des alg\`ebres de Lie
r\'esolubles exponentielles. Dans le cas nilpotent nous montrons
que toute orbite coadjointe est la vari\'et\'e de Poisson
caract\'eristique d'un module topologique bien choisi.
\par}
\eject
\null
\titre{Table des mati\`eres}
\alinea{Introduction\hfill\rm 2}
\alinea{I. Objets g\'eom\'etriques li\'es aux alg\`ebres d\'eform\'ees\hfill\rm 6}
I.1. Modules topologiquement libres \hfill 6

I.2. Id\'eaux divisibles \hfill 7

I.3. Annulateurs \hfill 8

I.4. Involutivit\'e \hfill 9

I.5. Vari\'et\'es caract\'eristiques \hfill 10

I.6. Involutions et $*$-repr\'esentations \hfill 12

I.7. Modules topologiquement libres convergents \hfill 14
\alinea{II. Application au cas des vari\'et\'es de Poisson lin\'eaires
            \hfill\rm 15}
II.1. L'alg\`ebre $(\Cal A,*)$  \hfill 15

II.2. Modules topologiquement libres convergents sur $(\Cal A,*)$ \hfill 16

II.3. Unitarit\'e   \hfill 17
\alinea{III. Exemples \hfill\rm 18}
III.1. Le groupe de Heisenberg \hfill 19

III.2. L'alg\`ebre de Lie filiforme de pas $n$  \hfill 21

III.3 L'alg\`ebre de Lie du groupe affine de la droite \hfill 22

III.4 Un exemple r\'esoluble exponentiel de dimension 3 \hfill 23

III.5 L'alg\`ebre de Lie du groupe diamant \hfill 24

III.6 L'alg\`ebre de Lie du groupe des d\'eplacements du plan \hfill 27

III.7 Le cas semi-simple : modules de Verma \hfill 30

III.8 Modules de Verma et r\'ealit\'e \hfill 31
\alinea{IV. Le cas r\'esoluble exponentiel \hfill\rm 33}
IV.1 Rappels sur la m\'ethode des orbites pour les groupes exponentiels
\hfill 33

IV.2 Construction d'une repr\'esentation $\pi_\nu$ de $\Cal A$  \hfill 36

IV.3 Calcul de la vari\'et\'e caract\'eristique \hfill 38

IV.4 Le cas nilpotent \hfill 39
\alinea{R\'ef\'erences\hfill\rm 39}
 \paragraphe{Introduction}
 \qquad
Une sous-vari\'et\'e involutive (ou co-isotrope) d'une vari\'et\'e de Poisson $V$ est une sous-vari\'et\'e plong\'ee $W$ telle que l'id\'eal des fonctions qui s'annulent sur $W$  est une sous-alg\`ebre de Poisson de $C^\infty(V)$. Dans le contexte de la quantification par d\'eformation, plusieurs auteurs (\cite {BGHHW}, \cite{CF2}) ont r\'ecemment propos\'e des m\'ethodes pour associer \`a une sous-vari\'et\'e involutive un id\'eal \`a gauche de l'alg\`ebre d\'eform\'ee $(C^\infty(V)[[\nu]],*)$, o\`u l'\'etoile-produit sur $V$ provient des constructions de M. Kontsevich \cite {K} ou de D. Tamarkin \cite{T}.
\ssq
Nous proposons dans cet article une d\'emarche inverse : d\'ecrire un proc\'ed\'e pour {\sl
d\'equantifier} certains modules sur une alg\`ebre d\'eform\'ee,
c'est-\`a-dire associer \`a un tel module une sous-vari\'et\'e
involutive et une sous-vari\'et\'e de Poisson de la vari\'et\'e de
Poisson sous-jacente. Le cadre $C ^\infty$ est mal adapt\'e, mais
le cadre analytique ou alg\'ebrique convient~: nous travaillons
d'abord sur le corps des complexes, puis, consid\'erant une involution naturelle sur l'alg\`ebre d\'eform\'ee nous mettons en \'evidence une classe de modules pour lesquels les objets ainsi mis en \'evidence vivent sur le corps des r\'eels. Ce sont les modules {\sl fortement pseudo-unitaires}, c'est-\`a-dire les modules munis d'une forme bilin\'eaire hermitienne non d\'eg\'en\'er\'ee compatible avec l'involution, \`a valeurs dans $\C[[\nu]]$ (o\`u $\nu$ est le param\`etre de d\'eformation), et telle que la forme bilin\'eaire quotient \`a valeurs dans $\C$ obtenue en $\nu=0$ est encore non-d\'eg\'en\'er\'ee.
\medskip
 \qquad
Soit $(V,\{,\})$ une vari\'et\'e de Poisson analytique r\'eelle
(resp. alg\'ebrique) c'est-\`a-dire une vari\'et\'e analytique
r\'eelle (resp. alg\'ebrique) munie d'un $2$-tenseur $P$ \`a
coefficients analytiques (resp. r\'eguliers) tel que le crochet
de Schouten $[P,P]$ s'annule. Le crochet de Poisson munit le
faisceau structural $\Cal O$ des germes de fonctions analytiques
(resp. r\'eguli\`eres) d'une structure de faisceau d'alg\`ebres de
Poisson. \ssq Nous nous limiterons au cas plat $V=\R^d$. En
complexifiant nous obtenons donc une structure de vari\'et\'e de
Poisson analytique complexe (resp. alg\'ebrique) sur $V^\C=\C^d$.
Nous noterons $A$ l'alg\`ebre de Poisson $\Cal O_V$ des fonctions
analytiques (resp. polynomiales) sur $\C^d$ (c'est l'espace des
sections globales du faisceau structural). M. Kontsevich \cite K a
construit un \'etoile-produit $\#$ sur $V$~~:
$$f\#g=\sum_{k\ge 0}\nu^kC_k(f,g),$$
o\`u les coefficients $C_k$ sont des
op\'erateurs bidiff\'erentiels d\'ecrits \`a l'aide de formules
compl\`e\-te\-ment explicites ne faisant intervenir que les
d\'eriv\'ees partielles des constantes de structure du $2$-tenseur
de Poisson (voir aussi \cite {AMM} et \cite {MT}). En particulier
 si le $2$-tenseur de Poisson est \`a coefficients analytiques
(resp. polynomiaux) alors l'\'etoile-produit est aussi \`a
coefficients analytiques (resp. polynomiaux). Autrement dit
l'\'etoile-produit $\#$ munit $\Cal A=A[[\nu]]$ d'une structure
d'alg\`ebre associative topologiquement libre sur $\C[[\nu]]$.
Cette alg\`ebre est naturellement filtr\'ee par $\Cal
A_n=\nu^n\Cal A$, et son gradu\'e associ\'e $\mop{Gr}\Cal A$ est
naturellement isomorphe \`a l'alg\`ebre de polyn\^omes $A[\nu]$ munie
du produit commutatif
de $A$ \'etendu par $\C[\nu]$-lin\'earit\'e. Enfin cet
\'etoile-produit est \`a coefficients r\'eels~: si $f$ et $g$ sont
des fonctions sur $V$ \`a valeurs r\'eelles, alors $f\#g$ est
aussi \`a valeurs r\'eelles. \ssq \ssq L'\'etoile-produit $\#$ de
M. Kontsevich est \'equivalent \`a un autre \' etoile-produit $*$
\cite {K \S\ 8}, \cite{MT}, \cite {CFT} ayant les m\^emes
propri\'et\'es, et v\'erifiant de plus la propri\'et\'e suivante~:
pour tout $f,g$ dans le centre de $(\Cal A,*)$ on a~:
$$f*g=fg.$$
Les deux alg\`ebres $(\Cal A,\#)$ et $(\Cal A,*)$ sont
bien entendu isomorphes. Nous appellerons le produit $*$ {\sl
\'etoile-produit de Duflo-Kontsevich} \cite {ABM}.
\ssq
Nous nous placerons exclusivement dans le cadre alg\'ebrique. Nous
d\'efinissons dans le premier paragraphe la vari\'et\'e caract\'eristique
$V(\Cal M)\subset V^\C$ d'un $\Cal A$-module topologiquement libre
$\Cal M$ ainsi que sa vari\'et\'e de Poisson caract\'eristique
$V\!A(\Cal M)$ en adaptant les d\'efinitions de \cite J de la fa\c con suivante~: consid\'erant l'annulateur $\mop{Ann} \Cal M$ du $\Cal A$-module $\Cal M$ nous d\'efinissons $V(\Cal M)$ comme l'ensemble des z\'eros communs de l'annulateur du $A$-module $M=\Cal M/\nu\Cal M$, et $V\!A(\Cal M)$ comme l'ensemble des z\'eros communs de $\mop{Ann}\Cal M/(\mop{Ann}\Cal M\cap\nu \Cal A)$. Ces deux
objets sont des sous-vari\'et\'es affines de $\C^n$ (i.e. d\'efinies par
l'annulation d'un nombre fini de
polyn\^omes).
\ssq
Nous montrons (\`a l'aide de \cite G) que
$V(\Cal M)$ est une sous-vari\' et\'e involutive de $V^\C$, que
$V\!A(\Cal M)$ est une sous-vari\'et\'e de Poisson de $V^\C$ (ce
qui justifie l'appellation), et que l'on a toujours l'inclusion~~:
$$V(\Cal M)\subset V\!A(\Cal M).$$ Apr\`es avoir montr\'e (\`a
l'aide de \cite {CF1}) que l'involution $f\longmapsto f^*$
d\'efinie
 par~:
$$f^*(\xi)=\overline{f(\overline \xi)}$$
est un
 anti-automorphisme de l'alg\`ebre $(\Cal A,*)$, nous introduisons
la notion de $\Cal A$-module {\sl fortement pseudo-unitaire}, et
nous montrons que la vari\'et\'e caract\'eristique d'un
$\Cal A$-module fortement pseudo-unitaire est
d\'efinie sur le corps des r\'eels, ainsi que sa vari\'et\'e de
Poisson caract\'eristique.
\ssq
Plus pr\'ecis\'ement, on prolonge d'abord la conjugaison complexe en un automorphisme de $\C[[\nu]]$ en d\'ecr\'etant que $\nu=i\hbar$ est imaginaire pur, c'est-\`a-dire $\overline\nu=-\nu$. Nous dirons que le module $\Cal M$ est pseudo-unitaire (ou, de fa\c con \'equivalente, que la repr\'esentation qui lui est associ\'ee est une $*$-repr\'esentation \cite {Bu-W}), s'il existe une forme sesquilin\'eaire non-d\'eg\'en\'er\'ee $<-,->_\nu$ \`a valeurs dans $\C[[\nu]]$ sur $\Cal M$ qui soit hermitienne, i.e.
$$<m,n>_\nu=\overline{<n,m>_\nu},\hskip 8mm m,n\in\Cal M$$
et compatible avec l'involution, i.e. v\'erifiant pour tout $a\in\Cal A$~:
$$<am,n>_\nu=<m,a^*n>_\nu.$$
La forme $<-,->_\nu$ induit par passage au quotient une forme hermitienne $<-,->_0$ sur $M$ \`a valeurs dans $\C$. Si cette forme est non-d\'eg\'en\'er\'ee nous dirons que le module $\Cal M$ est fortement pseudo-unitaire (ou, de fa\c con \'equivalente, que la repr\'esentation associ\'ee est une $*$-repr\'esentation fortement non-d\'eg\'en\'er\'ee). Contrairement à \cite {Bu-W}, nous ne faisons pas forc\'ement d'hypoth\`ese de positivit\'e sur la forme sesquilin\'eaire. Un module fortement pseudo-unitaire muni d'une forme d\'efinie positive sera dit {\sl fortement unitaire}.
\msq
Dans le deuxi\`eme paragraphe nous adaptons ce cadre aux vari\'et\'es de Poisson lin\'eaires. Soit donc $\g g$ une alg\`ebre de Lie r\'eelle de dimension finie, et soit $V=\g g^*$. L'alg\`ebre d\'eform\'ee $\Cal A=S(\g g)[[\nu]]$ s'identifie, via la version formelle de l'isomorphisme de Duflo, \`a l'alg\`ebre enveloppante formelle complexifi\'ee~:
$$\Cal U_\nu(\g g^\C)=T(\g g^\C)[[\nu]]/<x\otimes y-y\otimes x-\nu [x,y]>.$$
Il est possible de sp\'ecialiser l'ind\'etermin\'ee $\nu$ en une valeur non nulle dans l'\'ecriture de l'\'etoile-produit \cite{K\S\ 8}~: introduisons la famille \`a un param\`etre complexe d'alg\`ebres de Lie $\g g^\C_{\nu_0}$, de m\^ eme espace sous-jacent $\g g^\C$ avec le crochet de Lie~:
$$[X,Y]_{\nu_0}=\nu_0[X,Y].$$
L'\'evaluation en $\nu=\nu_0$ fournit une loi non-commutative $*_{\nu_0}$ sur $S(\g g^\C)$, qui est la multiplication de l'alg\`ebre enveloppante de $\g g_{\nu_0}^\C$ transport\'ee par l'isomorphisme de Duflo. Nous introduisons au paragraphe I.7 la notion de {\sl module topologiquement libre faiblement convergent\/}, qui nous permet de sp\'ecialiser l'ind\'etermin\'ee $\nu$ en une valeur non nulle \'egalement au niveau des repr\'esentations~: de fa\c con pr\'ecise un module topologiquement libre faiblement convergent sur $\Cal A$ est un $\Cal A$-module topologiquement libre $\Cal M=M[[\nu]]$ sur $\Cal A$, o\`u $M$ est un espace topologique localement convexe s\'epar\'e, tel qu'il existe $R>0$ tel que pour tout $a\in \Cal A_0$ et pour tout $m\in M$ la s\'erie enti\`ere $\pi_\nu(a)m$ est faiblement convergente de rayon $R$. Ici $\Cal A_0$ d\'esigne la sous-$\C[\nu]$-alg\`ebre de $\Cal A$ engendr\'ee par $A$ et $\pi_\nu$ la repr\'esentation associ\'ee.
\ssq
Soit maintenant $(\g g_\hbar)_{\hbar\in\R}$ la famille \`a un param\`etre r\'eel d'alg\`ebres de Lie r\'eelles de m\^eme espace sous-jacent $\g g$, avec le crochet~:
$$[X,Y]_\hbar=\hbar[X,Y].$$
Nous montrons (Proposition II.2.1 et Th\'eor\`eme II.3.1) l'\'equivalence entre la notion de $\Cal A$-module fortement unitaire faiblement convergent de rayon $R$ et la notion de famille \`a un param\`etre $(\rho_\hbar)_{\hbar\in ]-R,R[}$ de repr\'esentations unitaires de $\g g_\hbar$ telle que la famille de repr\'esentations de $\Cal U(\g g_\hbar)\simeq \bigl(S(\g g),*_{\hbar}\bigr)$ associ\'ee d\'epende faiblement analytiquement du param\`etre $\hbar$. Le passage d'une notion \`a l'autre est donn\'e par la formule~:
$$\rho_\hbar(X)=-i\pi_\nu(X)\restr{\nu=i\hbar}$$
pour tout $X\in\g g$. Le facteur $i$ s'explique par le fait que $\hbar$ est r\'eel alors que l'ind\'etermin\'ee $\nu$ est formellement imaginaire pure.
\msq
Le troisi\`eme paragraphe est consacr\'e \`a quelques exemples dans le cadre Poisson-lin\'eaire, plus pr\'ecis\'ement dans les cadres nilpotent et r\'esoluble, mis \`a part les modules de Verma trait\'es en conclusion.
\msq
Dans la derni\`ere partie, \`a l'aide de la m\'ethode des
orbites de Kirillov \cite {Ki} et des r\'esultats de N.V. Pedersen \cite {Pe1}, \cite{Pe2}  nous d\'eterminons, dans le cas des alg\`ebres de Lie r\'esolubles exponentielles la vari\'et\'e caract\'eristique d'un module fortement unitaire obtenu par induction unitaire d'une polarisation r\'eelle quelconque.
Nous montrons enfin que lorsque $V$ est le dual d'une
alg\`ebre de Lie nilpotente toute feuille symplectique (c'est-\`
a-dire toute orbite coadjointe) peut se voir comme la vari\'et\'e
de Poisson caract\'eristique d'un $\Cal A$-module fortement unitaire $\Cal
M$ bien choisi.
\ssq
Cette \'etude montre
que le cadre trac\'e au \S\ I pour une vaste classe de
vari\'et\'es de Poisson s'accorde bien avec la m\'ethode des
orbites de Kirillov dans le cas des vari\'et\'es de Poisson
lin\'eaires.
\smallskip
{\sl Remarque\/}~: contrairement \`a l'anneau des polyn\^omes, l'anneau des fonctions holomorphes sur $\C^n$ n'est pas noeth\'erien. En revanche l'anneau des germes de fonctions holomorphes en un point donn\'e est noeth\'erien. Dans l'optique d'une adaptation de nos m\'ethodes au cadre analytique il serait donc int\'eressant de ``faisceautiser'' la construction pr\'esent\'ee ici, c'est-\`a-dire de d\'equantifier des faisceaux de modules topologiquement libres sur le faisceau $(\Cal O[[\nu]],*)$ d'alg\`ebres d\'eform\'ees (qui a un sens puisque l'étoile-produit est bidiff\'erentiel). Cette m\^eme approche devrait permettre d'adapter le pr\'esent travail \`a des vari\'et\'es alg\'ebriques lisses plus g\'en\'erales \cite {Y}.
\smallskip
{\bf Remerciements}~: Nous remercions Georges Pinczon d'avoir attir\'e notre attention sur le point ci-dessus, et Charles Torossian pour ses remarques pertinentes.

\vskip 5mm
 \paragraphe{I. Objets g\'eom\'etriques associ\'es aux alg\`ebres
d\'eform\'ees}

\qquad Nous gardons les notations de l'introduction en nous limitant exclusivement au cadre alg\'ebrique. Nous introduisons
la notion d'id\'eal divisible de l'alg\`ebre d\'eform\'ee $\Cal A$, puis
nous introduisons, en nous inspirant de \cite J, la notion de vari\'et\'e
caract\'eristique pour un $\Cal A$-module topologiquement libre. C'est
un ferm\'e de Zariski de $V$. Nos vari\'et\'es caract\'eristiques ne sont pas forc\'ement
coniques (contrairement \`a ce qui se passe dans \cite J ou \cite {GM}),
ceci gr\^ace \`a l'introduction du param\`etre de d\'eformation dans la
construction. La notion de {\sl pseudo-unitarit\'e forte} permet de faire vivre ces vari\'et\'es caract\'eristiques sur le corps des r\'eels.
\alinea{I.1. Modules topologiques sur l'anneau des s\'eries formelles}
\qquad Nous reprenons les d\'efinitions de \cite {Ka \S \ XVI}, \cite {CFT \S\ A.1} et \cite {EK \S\ 2.1}. Soit $k[[\nu]]$ l'anneau des séries formelles sur un corps quelconque $k$. On munit cet anneau de la topologie $\nu$-adique, d\'efinie par la distance ultram\'etrique~:
$$d(a,b)=2^{-\smop{val}(a-b)},$$
avec $\mop{val}a=\mop{sup}\{j,a\in\nu^jk[[\nu]]\}$. Cette distance fait de $k[[\nu]]$ un anneau topologique complet. Sur tout $k[[\nu]]$-module $\Cal M$ on met une topologie invariante par translation en d\'ecidant que les $\nu^j\Cal M,j\in\N$ forment une base de voisinages de z\'ero. Cette topologie est s\'epar\'ee si et seulement si l'intersection des $\nu^j\Cal M$ est r\'eduite \`a $\{0\}$. Dans ce cas on peut d\'efinir la valuation~:
$$\mop{val}m=\mop{sup}\{j,m\in\nu^j\Cal M\}$$
et d\'efinir la topologie par la distance ultram\'etrique associ\'ee $d(m,m')=2^{-\smop{val}(m-m')}$.
\ssq
On dira qu'un $k[[\nu]]$-module $\Cal M$ est {\sl sans torsion\/} si l'action de $\nu$ est une injection de $\Cal M$ dans $\Cal M$. On rappelle \'egalement \cite {EK\S\ 2.1} qu'un $k[[\nu]]$-module {\sl topologiquement libre\/} est un module $\Cal M$ isomorphe \`a $M[[\nu]]$ pour un certain espace vectoriel $M$.
\prop{I.1.1 (\rm \cite{Ka} proposition XVI.2.4 et \cite{CFT} lemma A1)}
Un $k[[\nu]]$-module $\Cal M$ est topologiquement libre si et seulement s'il est s\'epar\'e, complet et sans torsion.
\ndem
Voir \'egalement \cite{B \S\ III.5.4}. Remarquons que pour un module topologiquement libre $\Cal M=M[[\nu]]$, l'espace vectoriel $M$ peut se voir \'egalement comme le quotient $\Cal M/\nu \Cal M$.
\ssq
Soit $\Cal M$ un $k[[\nu]]$-module, et soit $\Cal N$ un sous-$k[[\nu]]$-module de $\Cal M$. Consid\'erons les $k$-espaces vectoriels $M=\Cal M/\nu\Cal M$ et $N=\Cal N/\nu\Cal N$. L'inclusion $i:\Cal N\inj 5\Cal M$ induit une application $k$-lin\'eaire~:
$$i_0:N\longrightarrow M.$$
On dira \cite {E \S\ 4.5} que le sous-$k[[\nu]]$-module $\Cal N$ est {\sl divisible\/} si l'application $i_0$ est injective. Autrement dit $\Cal N$ est divisible si et seulement si $\nu \Cal N=\Cal N\cap\nu\Cal M$.
Un exemple simple de sous-$k[[\nu]]$-module non divisible est $\nu M[[\nu]]\subset M[[\nu]]$. Pour tout $k[[\nu]]$-module $\Cal M$, on \'ecrira aussi ``$m=O(\nu^j)$ dans $\Cal M$'' pour $m\in\nu^j\Cal M$. Un sous-$k[[\nu]]$-module $\Cal N$ de $\Cal M$ est donc divisible si et seulement si pour tout $m\in\Cal N$, $m=O(\nu)$ dans $\Cal M$ implique $m=O(\nu)$ dans $\Cal N$.
\ssq
L'alg\`ebre associative $\Cal A$ d\'efinie dans
l'introduction est par construction un $\Cal \C[[\nu]]$-module
topologiquement libre. Nous appellerons {\sl $\Cal A$-module
topologique\/} un $\Cal \C[[\nu]]$-module $\Cal M$ muni d'une application $\Cal
\C[[\nu]]$-bilin\'eaire~:
$$\eqalign{\Phi :\Cal A\times\Cal M &\longrightarrow \Cal M \cr (\varphi,m)
&\longmapsto \pi_\nu(\varphi)m\cr}$$
faisant de $\Cal M$ un $\Cal A$-module.
\prop{I.1.2}
Soit $\Cal M$ un $\Cal A$-module topologique. Alors l'application $\C[[\nu]]$-bilin\'eaire $\Phi : \Cal A\times\Cal M \rightarrow \Cal M$ d\'efinissant le module est continue pour les topologies $\nu$-adiques de $\Cal A$ et de $\Cal M$.
\dem
Soient $m\in\Cal M$ et $a\in A$. la famille $\Cal W_j=\Phi(a,m)+\nu^j\Cal M$ forme une base de voisinages de $\Phi(a,m)=\pi_\nu(a)m$. Consid\'erant les voisinages $\Cal U_j=a+\nu^j\Cal A$ et $\Cal V_j=m+\nu^j\Cal M$ de $a$ et de $m$ respectivement, il est alors clair que l'image par $\Phi$ du produit $\Cal U_i\times\Cal V_j$ est incluse dans $\Cal W_j$, ce qui montre la continuit\'e. En particulier la multiplication de $\Cal A\times\Cal A$ dans $\Cal A$ est $\C[[\nu]]$-bilin\'eaire continue pour la topologie $\nu$-adique, et fait donc de $\Cal A$ une alg\`ebre topologique.
\qed
 \ssq
Soient $\Cal M_1$ et $\Cal M_2$ deux
$\Cal A$-modules topologiques. Un morphisme de $\Cal A$-modules
topolo\-giques (ou {\sl op\'erateur d'entrelacement}) de $\Cal
M_1$ vers $\Cal M_2$ est une application $\C[[\nu]]$-lin\'eaire continue
commutant aux actions des \'el\'ements de $\Cal A$. On dit que $\Cal
M_1$ et $\Cal M_2$ sont {\sl \'equivalents} s'il existe un op\'erateur
d'entrelacement bijectif bicontinu de $\Cal M_1$ sur $\Cal M_2$.
\goodbreak
\alinea{I.2. Id\'eaux divisibles}
\qquad Comme $\Cal A=A[[\nu]]$ est topologiquement libre on identifiera $A$ avec $\Cal A/\nu\Cal A$. Soit $\Cal J$ un id\'eal \`a gauche de $\Cal A$. Il est
imm\'ediat de voir que $J=\Cal J/(\Cal J\cap\nu\Cal A)$ est un id\'eal de
l'alg\`ebre commutative $A$. La situation est bien s\^ur la
m\^eme avec les id\'eaux \`a droite.
\prop{I.2.1}
Soit $\Cal J$ un id\'eal bilat\`ere divisible de $\Cal A$. Alors $J=\Cal
J/(\Cal J\cap\nu\Cal A)$ est un id\'eal de Poisson de $A$.
\dem
Soient $a_0\in J$ et $b_0\in A$. Soit $a=a_0+\nu
a_1+\cdots\in \Cal J$ repr\'esentant $a_0$ et soit $b=b_0+\nu b_1+\cdots\in \Cal A$ repr\'esentant $b_0$. Par divisibilit\'e de $\Cal J $ on a~:
$${1\over
\nu}(a*b-b*a)\in\Cal J,$$
d'o\`u on d\'eduit, en consid\'erant le
terme constant, que $\{a_0,b_0\}$ appartient \`a $J$. \qed
Remarquons que, comme $\Cal J$ est divisible, on peut aussi identifier $J$ \`a $\Cal J/\nu\Cal J$.
\alinea{I.3. Annulateurs}
Soit $\Cal M$ un $\Cal A$-module topologique. On d\'efinit l'{\sl annulateur\/} $\mop{Ann} \Cal M$ du module $\Cal
M$ comme l'ensemble des $\varphi\in\Cal A$ tels que
$\pi_\nu(\varphi)m=0$ pour tout $m\in\Cal M$. On voit
imm\'ediatement que alors
$\mop{Ann} \Cal M$ est un id\'eal bilat\`ere de $\Cal A$, et que $\mop{Ann} \Cal M$ est divisible si $\Cal M$ est sans torsion.
\prop{I.3.1}
Soit $\Cal M$ un $\Cal A$-module. Alors $M=\Cal M/ \nu
\Cal M$ est un module sur l'alg\`ebre commutative $A=\Cal A/\nu \Cal
A$.
\dem C'est imm\'ediat. On notera $\pi_0$ la
repr\'esentation de $A$ dans le module $M$ associ\'ee. \qed
\prop{I.3.2}
Soit $\Cal M$ un $\Cal A$-module topologique. Alors
l'annulateur du $A$-module $M=\Cal M/\nu\Cal M$ est stable par le crochet de Poisson de
$A$.
\dem On note $\pi_\nu$ la repr\'esentation de $\Cal A$ dans le module $\Cal M$. L'annulateur de $M$ peut se voir comme l'ensemble des $f\in
A$ tels que $ \pi_\nu(f)u=O(\nu)$ pour tout $u\in \Cal M$. Pour tout
$f,g\in \mop{Ann } M$ on a dans $\Cal M$~:
$$\pi_\nu(f*g)u=\pi_\nu(f)\pi_\nu(g)u=O(\nu^2),$$
ce qui nous donne~:
$$\eqalign{\pi_\nu(\{f,g\})m &={1\over\nu}\pi_\nu(f*g-g*f).m+O(\nu)
\cr &=O(\nu),\cr}$$
d'o\`u le r\'esultat. \qed
{\sl Remarque}~: L'anneau
$A$ \'etant noeth\'erien \cite {B}, l'annulateur de $M$ est finiment engendr\'e.
\alinea{I.4. Involutivit\'e}
\qquad Soit $V=\R^d$ une vari\'et\'e de Poisson plate alg\'ebrique r\'eelle (comme dans l'introduction), et soit
$W\subset V^\C$ une sous-vari\'et\'e affine. Soit $I(W)$
l'id\'eal de $A$ constitu\'e par les fonctions qui s'annulent sur
$W$. On dira que $W$ est {\sl involutive\/} ou {\sl co-isotrope\/} si l'id\'eal $I(W)$ est stable par le crochet de Poisson. C'est \'equivalent \`a l'annulation du $2$-tenseur de Poisson sur la deuxi\`eme puissance ext\'erieure du fibr\'e conormal de $W$ \cite {BGHHW}, \cite{CF2}.
\prop{I.4.1}
Soit $W$ une sous-vari\'et\'e affine de $V^\C$. Soit $W^r$ la partie non-singuli\`ere de $W$. Si $W$
est involutive, alors pour toute feuille symplectique $S$ de
$V^\C$ coupant $W^r$ transversalement l'intersection $W^r\cap S$
est une sous-vari\'et\'e co-isotrope de $S$, c'est-\`a-dire que
pour tout $x\in W^r\cap S$ on a~~: $$(T_x(W^r\cap
S))^\omega\subset T_x(W^r\cap S),\eqno{(\I.4.1)}$$ o\`u l'exposant
$\omega$ d\'esigne l'orthogonal dans $T_xS$ pour la forme
symplectique. \dem Nous adaptons la d\'emonstration de la
proposition 19 de \cite {GM \S\ II I.4}. Soit $x\in W^r\cap S$.
Soit $P_x$ le 2-tenseur de Poisson en $x$, et soit $\wt
P_x~:T^*_xV^\C\to T_xV^\C$ l'application lin\'eaire
antisym\'etrique associ\'ee, d\'efinie par~:
$$<\eta,\,\widetilde P_x(\xi)>=P_x(\xi,\,\eta).\eqno{(\I.4.2)}$$
L'image de $\wt
P_x$ est pr\'ecis\'ement $T_xS$.

\lemme{I.4.2}

On a l'\'egalit\'e~:
$$T_x(W^r\cap S)^\omega=\wt
P_x\bigl(T_x(W^r\cap S)^\perp\bigr),\eqno{( \I.4.3)}$$
o\`u l'exposant
$\perp$ d\'esigne l'orthogonal d'un sous-espace dans le dual.
\dem Soit $\xi\in T_x(W^r\cap S)^\perp$. Alors pour tout $Y\in
T_x(W^r\cap S)$ on a~: $$\omega(\wt P_x(\xi),\,Y)=<\xi,\,Y>=0,
\eqno{(\I.4.4)}$$ ce qui montre l'inclusion~: $$\wt
P_x\bigl(T_x(W^r\cap S)^\perp\bigr)\subset T_x(W^r\cap
S)^\omega.$$ Pour montrer l'inclusion inverse on consid\`ere un
$X$ dans $T_x(W^r\cap S)^\omega$. C'est l'image par $\wt P_x$ d'un
\'el\'ement $\xi$ de $T_x^*V $, et il est imm\'ediat d'apr\`es
(I.4.4) que $\xi$ appartient \`a $T_x(W ^r\cap S)^\perp$.
\qed
{\sl Fin de la d\'emonstration de la proposition I.4.1}~:
Soit $X\in T_xS$. Par d\'efinition d'une
feuille symplectique il existe $\varphi\in A$ telle que $X$ co\"\i
ncide avec le champ hamiltonien $H_\varphi(x)$. Alors $X$
appartient \`a $( T_x(W^r\cap S))^\omega$ si et seulement si pour
tout $Y\in T_x(W^r\cap S) $ on a~:
$$\omega(X,Y)=Y.\varphi(x)=d\varphi(x)(Y)=0.$$
\ssq
Maintenant, l'espace $B_x$ des $d\varphi (x)$ o\`u
$\varphi\in I(W)$ est l'orthogonal de $T_xW^r$. Gr\^ace \`a la
condition de transversalit\'e nous pouvons donc \'ecrire~:
$$T_x(W^r\cap S)^\perp=(T_xW^r\cap
T_xS)^\perp=B_x+(T_xS)^\perp,\eqno {(\I.4.5)}$$
Soient donc $X,Y\in T_x(W^r\cap S)^\omega$. D'apr\`es (I.4.5) et le
lemme I.4.2 il existe $\varphi,\psi\in I(W)$ telles que $X=\wt
P_x(d\varphi (x))=H_ \varphi(x)$ et $Y=H_ \psi(x)$. On a alors gr\^ace \`a l'involutivit\'e~:
$$\omega(X,Y)=\{\varphi,\psi\}(x)=0,$$
Ce qui termine la d\'emonstration.
 \qed
{\sl Remarque}~: on a l'\'egalit\'e~:
$$T_x(W^r\cap S)^\omega=\wt
P_x\bigl((T_xW^r)^\perp\bigr).\eqno{(\I.4.7)}$$
On a une r\'eciproque partielle \`a la
proposition I.4.1~:

\prop{I.4.3}

Soit $W$ une sous-vari\'et\'e affine de $V^\C$. Soit $W^r$ la partie non-singuli\`ere de $W$.
Supposons qu'il existe un ensemble Zariski-dense $\Cal U$ de $W^r$ tel que~:

1) pour tout $x\in\Cal U$, l'intersection de $W^r$ avec la
feuille symplectique $S_x$ passant par $x$ est transverse.

2) $W^r\cap S_x$ est co-isotrope dans $S_x$.
\medskip
Alors $W$ est involutive.

\dem On voit facilement que sous
les hypoth\`eses de la proposition, si $f$ et $g$ appartiennent \`a $I(W)$
alors $\{f,g\}$ s'annule sur l'ensemble Zariski-dense $\Cal U$ de $W^r$, et donc sur $W$
tout entier. Rappelons \cite{V} que m\^eme dans le cas alg\' ebrique que nous
consid\'erons, les feuilles symplectiques ne sont pas en g\'en\'eral des
sous-vari\'et\'es alg\'ebriques de $V^\C$.
Un exemple de cette situation est donn\'e au paragraphe III.4.\qed

\alinea{I.5. Vari\'et\'es caract\'eristiques}

\qquad Soit $\Cal M$ un $\Cal A$-module topologique. La
{\sl vari\'et\'e caract\'eristique\/} $V(\Cal M)$ de $\Cal M$ est
d\'efinie comme l'ensemble des z\'eros communs de l'annulateur du
$A$-module $M=\Cal M/ \nu\Cal M$. Si $\Cal M$ est sans torsion, on
appellera {\sl vari\'et\'e de
Poisson caract\'eristique de $\Cal M$} et, comme dans \cite J, on
notera $V\! A(\Cal M)$ l'ensemble des z\'eros communs de l'id\'eal
$\mop{Ann}\Cal M/(\mop{Ann}\Cal M\cap \nu\Cal A)$ de $A$. Cette terminologie
est justifi\'ee au corollaire I.5.3 ci-dessous.
\ssq
Comme $A$ est
un anneau commutatif noeth\'erien, le gradu\'e associ\'e $A[\nu]$
de $\Cal A$ est aussi noeth\'erien \cite
{B \S\ III.2 corollaire 1}. On peut donc appliquer le th\'eor\`eme
d'int\'egrabilit\'e des vari\'et\'es caract\'eristiques de Gabber
\cite {G, Theorem I}~:
\th{I.5.1 (\rm Int\'egrabilit\'e des caract\'eristiques ~: O.
Gabber)}
Supposons que $\Cal M$ soit un $\Cal A$-module finiment
engendr\'e, et soit $M=\Cal M/\nu\Cal M$. Alors le radical $J(\Cal
M)$ de $\mop{Ann} M$ est stable par le crochet de Poisson.
\ndem
{\sl Remarque}~: $J(\Cal M)$ est l'id\'eal des \'el\'ements de
$A$ qui s' annulent sur la vari\'et\'e caract\'eristique $V(\Cal
M)$. Le th\'eor\`eme dit donc que $V(\Cal M)$ est involutive.
\th{I.5.2}
L'ensemble des z\'eros communs d'un id\'eal de Poisson est une
sous-vari\'et\'e de Poisson de $V^\C$.
\dem
Soit $J$ un id\'eal de Poisson de $A$, et soit
$f\in J$. Soit $x$ annulant tous les \'el\'ements de $J$, et soit
$y=\phi_t(x)$, o\`u $(\phi_t)_{|t|<\varepsilon}$ est le flot d'un
champ de vecteurs hamiltonien $H_g, g\in A$. Par analyticit\'e du
flot, on a pour $|t|$ assez petit~:
$$(f\circ \phi_t)(x)=\sum_{k=0}^\infty {1\over k!}(\mop{ad}^k g.f)(x)t^k,\eqno{( \I.5.1)}$$
o\`u le membre de droite est convergent. Or
le terme~:
$$(\mop{ad}^k g.f)(x)=\{g,\{g,\ldots \{g,f\}\ldots\}\}(x)$$
s'annule, puisque $J$ est un id\'eal de Poisson. Donc $(f\circ \phi_t)(x) $
s'annule pour $t$ assez petit. On en d\'eduit que la feuille symplectique
passant par $x$ est enti\`erement contenue dans l'ensemble des z\'eros
communs de $J$. \qed

\cor{I.5.3}
La vari\'et\'e de Poisson caract\'eristique $V\!A(\Cal M)$ d'un $\Cal A$-
module topologique sans torsion $\Cal M$ est une sous-vari\'et\'e de
Poisson de $V^\C$.
\dem
Comme $\Cal M$ est sans torsion, l'id\'eal $\mop{Ann}\Cal M$ de $\Cal A$ est
divisible, et donc d'apr\`es la proposition I.2.1 $\mop{Ann}\Cal
M/(\mop{Ann}\Cal M\cap\nu\Cal A)$ est un id\'eal de Poisson de $A$. \qed
\prop{I.5.4}
Pour tout $\Cal A$-module topologique sans torsion $\Cal M$ on a
l'inclusion~: $$V(\Cal M)\subset V\!A(\Cal M).\eqno{(\I.5.2)}$$
\dem Ce r\'esultat d\'ecoule naturellement des caract\'erisations
suivantes~:
$$\eqalign{\mop{Ann} M = \mop{Ann}(\Cal M/\nu\Cal
M)&=\{\varphi_0\in A/\pi_\nu(\varphi_0)=O(\nu)\} \cr
\mop{Ann}\cal M/(\mop{Ann}\cal M\cap \nu\Cal A) & =\{\varphi_0\in A/\exists
\varphi=\varphi_0+\nu\varphi_1+.. .\in \Cal
A/\pi_\nu(\varphi)=0\}.\cr}\eqno{(\I.5.3)}$$
Il est clair que le
deuxi\`eme id\'eal est contenu dans le premier, d'o\` u l'inclusion
inverse des vari\'et\'es caract\'eristiques. \qed
{\sl Remarque}~:
les annulateurs de deux modules topologiquement libres \'equivalents sont
\'egaux. La vari\'et\'e de Poisson caract\'eristique $V \!A(\Cal M)$ ne
d\'epend donc que de la classe d'\'equivalence de $\cal M $. Il n'en va
pas de m\^eme de la vari\'et\'e caract\'eristique $V(\Cal M )$.
\alinea{I.6. Involution et $*$-repr\'esentations}
\qquad A. Cattaneo et G. Felder ont remarqu\'e \cite {CF1 \S\ 2} que
l'alg\`ebre $(\Cal A,\#)$ est naturellement munie d'une involution. Nous
adaptons ici leurs arguments, et nous allons montrer que cette
involution est \'egalement une involution pour l'\'etoile-produit~$*$ de
Duflo-Kontsevich. On \'etend la conjugaison complexe en un
automorphisme involutif de $\C[[\nu] ]$ en posant simplement
$\overline \nu=-\nu$. On consid\`ere donc $\nu$ comme ``imaginaire
pur'', et nous poserons \'egalement $\nu=i\hbar$, o\`u
$\hbar=\overline\hbar$ peut \^etre consid\'er\'e comme ``r\'eel''.

\prop{I.6.1}

L'involution semi-lin\'eaire $f\mapsto f^*$ de $\Cal A$ d\'efinie par~:
$$f^*(\xi)=\overline{f(\overline\xi)}\eqno{(\I.6.1)}$$
 est un anti-automorphisme de l'alg\`ebre $(\Cal A,\#)$.
\dem Dans l'expression de
l'\'etoile-produit on peut remplacer le param\`etre $\nu$ par n'importe
quelle s\'erie formelle sans terme constant. En particulier on peut
remplacer $\nu$ par $\hbar=-i\nu$ (qui v\'erifie alors $
\overline\hbar=\hbar$). L'\'etoile-produit $\#$ \'etant d\'efini par des
op\'erateurs bi-diff\'erentiels \`a coefficients r\'eels, on v\'erifie
facilement que l'on a~:
$$(f\#_\hbar g)^*=f^*\#_\hbar
g^*.\eqno{(\I.6.2)}$$
 d'o\`u l'on d\'eduit~:
$$(f\#_{\nu} g)^*=f^*\#_{-\nu}g^*.\eqno{(\I.6.3)}$$
Enfin, on sait (voir la remarque \`a la fin du
\S\ 2 de \cite{CF1}) que l' \'etoile-produit de Kontsevich v\'erifie la
propri\'et\'e de parit\'e altern\'ee de ses coefficients, c'est-\`a-dire~:
$$f\#_{-\nu}g=g\#_{\nu} f.\eqno{(\I.6.4)}$$
D'apr\`es (I.6.3) on a donc~:
$$(f\#_{\nu} g)^*=g^*\#_{\nu}f^*,\eqno{(\I.6.5)}$$
d'o\`u le r\'esultat.
\qed

\prop{I.6.2}

L'involution $f\mapsto f^*$ est \'egalement un anti-automorphisme
de l'alg\`ebre $(\Cal A,*)$ munie de l'\'etoile-produit de
Duflo-Kontsevich. Elle se restreint \`a $A$ et sa restriction est
un (anti-) automorphisme de l'alg\`ebre commutative $A$. \dem Soit
$D=I+\nu D_1+\cdots$ l'op\'erateur diff\'erentiel formel qui
r\'ealise l'\'equiva\-lence entre les deux \'etoile-produits~:
$$f*g=D(D\inver f\# D\inver g).\eqno{(\I.6.6)}$$
Cet op\'erateur d'\'equivalence $D$ est \`a coefficients r\'eels,
et commute donc avec l'involu\-tion $f\mapsto f^*$. L'involution
semi-lin\'eaire $f\mapsto f^0$ d\'efinie par~:
$$f^0=D[(D\inver f)^*]\eqno{(\I.6.7)}$$
co\"\i ncide donc avec $f\mapsto f^*$.
 Dans le cadre
des vari\'et\'es de Poisson lin\'eaires (voir \S\ II.1 ci-dessous) on a
m\^eme l'\'egalit\'e entre les deux \'etoile-produits. \ssq L'involution se
restreint de mani\`ere \'evidente \`a $A$. De plus l'involution respecte
$\nu\Cal A$, et l'involution ainsi d\'efinie sur $\Cal A/ \nu \Cal A$
correspond \`a cette restriction via l'isomorphisme canonique de $A$ sur
$\Cal A/\nu \Cal A$.
\qed
\msq
Soit $\pi_\nu$ une repr\'esentation de
$(\Cal A,*)$ dans un module topologique $\Cal M$. A la suite de \cite {Bu-W} on dit que $\pi_\nu$ est une {\sl $*$-repr\'esentation}
de $\Cal A$ dans $\Cal M$ si on a~:
$$<\pi_\nu(f)u,v>_\nu=<u,\pi_\nu(f)^* v>_\nu, \hskip 8mm u,v\in\Cal M,\ f\in\Cal A\eqno{(\I.6.8)}$$
 pour une certaine
forme sesquilin\'eaire $<-,->_\nu$ hermitienne non-d\'eg\'en\'er\'ee (sans hypoth\`ese de positivit\'e) sur $\Cal M$ \`a valeurs dans $\C[[\nu]]$. On dira que l'$*$-repr\'esentation est {\sl unitaire\/} si de plus cette forme sesquilin\'eaire est d\'efinie positive, c'est-\`a-dire $<x,x>_\nu>0$ pour tout $x\in\Cal M$ non nul, l'ordre sur $\R[[\nu]]$ \'etant l'ordre lexicographique.
\ssq
Cette forme sesquilin\'eaire d\'efinit par passage au quotient une forme sesquilin\'eaire sur $M=\Cal M/\nu\Cal M$ \`a valeurs dans $\C$. Nous supposerons que la forme $<-,->_\nu$ est {\sl fortement non-d\'eg\'en\'er\'ee}, c'est-\`a-dire que nous supposerons aussi non-d\'eg\'en\'er\'ee la forme quotient sur $M$. Dans ce cas on dira que l'$*$-repr\'esentation est fortement non d\'eg\'en\'er\'ee. Une repr\'esentation fortement non d\'eg\'en\'er\'ee et unitaire sera dite fortement unitaire.
\ssq
Un $\Cal A$-module unitaire (resp. pseudo-unitaire) sera par d\'efinition un module topologi\-que muni d'une repr\'esentation unitaire (resp. d'une $*$-repr\'esentation) de $\Cal A$. Un $\Cal A$-module fortement unitaire (resp. fortement pseudo-unitaire) sera par d\'efinition un module topologique muni d'une repr\'esentation fortement unitaire (resp. d'une $*$-repr\'esentation fortement non d\'eg\'en\'er\'ee) de $\Cal A$.
\prop{I.6.3}
Soit $\pi_\nu$ une $*$-repr\'esentation fortement non-d\'eg\'en\'er\'ee de $(\Cal A,*)$ dans un module
topologique $\Cal M$. Alors l'annulateur de $M=\Cal M/\nu
\Cal M$ est engendr\'e par un nombre fini d'\'el\'ements auto-adjoints
de $A$. Il en est de m\^eme de $\mop{Ann}\Cal M/(\mop{Ann}\Cal
M\cap \nu\Cal A)$.

\dem Tout \'el\'ement $f$ de $A$ s'\'ecrit de mani\`ere unique~:
$$f=f^++if^-,\eqno{(\I.6.9)}$$ o\`u $f^+$ et $f^-$ sont auto-adjoints. On a~:
$$f^+={1\over2}(f+f^*),\hskip 8mm f^-={1\over 2i}(f-f^*).\eqno{(\I.6.10)}$$
Gr\^ace \`a la non-d\'eg\'en\'erescence du produit scalaire sur $M$ on voit que
si $f$ appartient \`a $\mop{Ann}M$, $f^*$ appartient aussi \`a $\mop{Ann}M$. Soit
$\{f_1,\ldots ,f_k\}$ un syst\`eme de g\'en\'erateurs de
$\mop{Ann}M$. Il est alors clair que
$\mop{Ann}M$ est engendr\'e par $\{f^+_1,\ldots ,f^+_k, f^-_1,\ldots ,f^ -
_k\}$. Le m\^eme raisonnement s'applique \`a $\mop{Ann}\Cal
M/(\mop{Ann}\Cal M\cap\nu\Cal A)$, comme on peut le voir en utilisant la
deuxi\`eme \'egalit\'e (I.5.3) et la pseudo-unitarit\'e forte de $\Cal M$. \qed
\cor{I.6.4}
Sous les hypoth\`eses ci-dessus, la vari\'et\'e caract\'eristique
$V(\Cal M)$ d'un module fortement pseudo-unitaire est r\' eelle, de m\^eme que sa vari\'et\'e de Poisson caract\'eristique $V\!A(\Cal M)$ lorsque $\Cal M$ est sans torsion.
\dem
La vari\'et\'e caract\'eristique $V(\Cal M)$ est d\'efinie par les
$2k$ \'equations~:
$$f^+_j( \xi)=f^-_j(\xi)=0,\ j=1,\cdots,
k.\eqno{(\I.6.11)}$$
Or les polyn\^omes $f^{\pm}_j$ sont bien \`a coefficients r\'eels par
d\' efinition de l'involution. Le raisonnement pour $V\!A(\Cal M)$
est identique. \qed
{\sl Remarque 1\/}~: La proposition I.6.3 et
le corollaire I.6.4 n'utilisent pas d'hypoth\`ese de positivit\'e sur le produit scalaire. Il faut noter que les notions de non-d\'eg\'en\'erescence introduites ici sont diff\'erentes de celles introduites dans \cite {Bu-W}, qui par ailleurs se placent d'embl\'ee dans le cas unitaire.
\medskip
{\sl Remarque 2}~: On notera indiff\'eremment $V(\Cal M)=V(\pi_\nu)$ pour la vari\'et\'e caract\'eristique, et $V\!A(\Cal M)=V\!A(\pi_\nu)$ pour la vari\'et\'e de Poisson caract\'eristique.
\alinea{I.7. Modules topologiquement libres convergents}
Nous aurons besoin par la suite de sp\'ecialiser le param\`etre de d\'eformation $\nu$
en une valeur complexe non nulle, en faisant converger les s\'eries formelles.
Soit $\Cal M=M[[\nu]]$ un module topologiquement libre sur l'alg\`ebre d\'eform\'ee
$\Cal A=A[[\nu]]$, et soit $\pi_\nu$ la repr\'esentation associ\'ee. On suppose de plus que $M$
 est un espace topologique localement convexe s\'epar\'e.
 Soit $\Cal A_0$ la sous-$\C[\nu]$-alg\`ebre de $\Cal A$ engendr\'ee
 par $A$. C'est l'ensemble des sommes~:
$$\sum_{j=0}^N\nu^j\alpha_j,$$
o\`u $n\in\N$ et chaque $\alpha_j$ est une somme de termes du type $a_1*\cdots *a_r$, avec $\uple ar\in A$. On dira
que $\Cal M=M[[\nu]]$ est faiblement convergent s'il existe $R>0$
tel que pour tout $a\in \Cal A_0$ et tout $m\in M$ la s\'erie
enti\`ere $\pi_\nu(a)m$ converge faiblement pour $\nu=\nu_0$ dans le disque de rayon $R$ vers un
vecteur de $M$ que l'on notera $\wt{\pi_{\nu_0}}(a)m$. L'unicit\'e de
cette limite faible est assur\'ee par le th\'eor\`eme de
Hahn-Banach. Le rayon de convergence $R_\Cal M$ du module est
alors d\'efini comme la borne sup\'erieure des rayons $R$
ci-dessus.
\prop{I.7.1} Un module topologiquement libre faiblement convergent
$\Cal M=M[[\nu]]$ de rayon $R_{\Cal M}$ induit une famille de
repr\'esentations $(\wt{\pi_{\nu_0}})_{\nu_0\in D(0,R_\Cal M)}$ de
$\Cal A_0$ dans $M$.
\dem Pour tout $a,b\in A$, pour tout $m\in M$ et pour tout $m'$
dans le dual topologique $M'$ on a \'egalit\'e entre s\'eries
formelles~:
$$<m',\,\pi_\nu(a*b)m>=<m',\,\pi_\nu(a)\pi_\nu(b)m>.\eqno{(\I.7.1)}$$
L'\'egalit\'e reste valide lorsque $a$ et $b$ sont dans $\Cal
A_0$. On a convergence de ces deux s\'eries enti\`eres en
$\nu=\nu_0\in D(0,R_{\Cal M})$ , et le membre de droite est aussi
la limite de la s\'erie enti\`ere
$<m',\,\pi_\nu(a)\wt{\pi_{\nu_0}}(b)m>$ en $\nu=\nu_0$. On a donc
pour tout $\nu_0\in D(0,R_{\Cal M})$ l'\'egalit\'e dans $M$~:
$$\wt{\pi_{\nu_0}}(a*b)m=\wt{\pi_{\nu_0}}(a)\wt{\pi_{\nu_0}}(b)m.\eqno{(\I.7.2)}$$
\qed
\paragraphe{II. Application au cas des vari\'et\'es de Poisson lin\'eaires}
\qquad Nous allons appliquer les r\'esultats pr\'ec\'edents aux
repr\'esentations des alg\`ebres de Lie. L'ingr\'edient essentiel
ici est la sp\'ecialisation du param\`etre de d\'eformation $\nu$
en un imaginaire pur quelconque, ce qui permet (Th\'eor\`eme
II.3.1) d'associer \`a un module topologiquement libre faiblement
convergent fortement unitaire sur l'alg\`ebre d\'eform\'ee une famille \`a un
param\`etre r\'eel $(\rho_\hbar)$ de repr\'esentations unitaires
d'alg\`ebres de Lie.
\ssq
Soit $V=\R^d$ une vari\'et\'e de Poisson
lin\'eaire. Alors le dual $V^*$ des formes lin\'eaires sur $V$
forme une sous-alg\`ebre de Lie $ \g g$ de l'alg\`ebre $A$ des
polyn\^omes sur $V$ munie du crochet de Poisson. On voit donc la
vari\'et\'e de Poisson $V$ comme le dual $\g g^*$ de l'alg\`ebre
de Lie $\g g$. Le crochet de Poisson de Kirillov-Kostant-Souriau
s'exprime pour $f,g\in A$ et $\xi \in \g g$ par la formule \cite
{W \S\ 3}~:
$$\{f,g\}(\xi)=<\xi,\ [df(\xi),\,dg(\xi)]>.$$
\alinea{II.1. L'alg\`ebre $(\Cal A,*)$}
Soit $\g g$ une alg\`ebre de Lie sur le corps $\R$, et soit $A=S(\g g^\C)$.
Il r\'esulte des travaux de B. Shoikhet sur l'annulation des poids associ\'es aux ``roues''
\cite {FS}, \cite S, que les deux \'etoile-produits $\#$ et
$*$ co\"\i ncident. L'alg\`ebre $(\Cal A,*)$ est isomorphe \`a "l'alg\`ebre enveloppante
formelle
complexifi\'ee"~: $$\Cal U_\nu(\g
g^\C)=T(\g g^{\C})[[\nu]]/<x\otimes y-y\otimes x-\nu[x,
y]>,\eqno{(\II.1.1)}$$ et on a pr\'ecis\'ement~:
$$f*g=\tau\inver(\tau f.\tau g).\eqno{(\II.1.2)}$$
Ici $\tau ~:\Cal A\to U_\nu(\g g^\C)$ est l'isomorphisme de
Duflo \cite {D1}~:
$$\tau=\sigma\circ J(D)^{1/2},\eqno{(\II.1.3)}$$
o\`u $\sigma$ est la
sym\'etrisation et $J(D)^{1/2}$ est l'op\'erateur diff\'erentiel d'ordre
infini \`a coefficients constants correspondant \`a la s\'erie formelle~:
$$J(x)^{1/2}=\Big (\mop{det}{\mop{sh}\ad {\frac\nu 2}x\over\ad
{\frac\nu 2}x}\Big )^{1/2}.\eqno{(\II.1.4)}$$

On notera $(S^n(\g g^\C))_{n\ge 0}$ la
filtration croissante usuelle de l'alg\`ebre sym\'etrique. On
peut sp\'ecialiser la valeur du param\`etre de d\'eformation
$\nu$~: en effet, pour tout $f,g$ dans $A $ la s\'erie en $\nu$
d\'efinissant $f*g$ est polynomiale en $\nu$ \cite {K \S\
8}, et donc peut s'\'evaluer en tout nombre complexe $\nu$~:
l'\'etoile-produit engendre donc une famille de
lois associatives non-commutatives $(*_\nu)$ sur $A$, le param\`etre $\nu$
parcourant l'ensemble des nombres complexes. Chacune de ces
alg\`ebres s'identifie via l'isomorphisme de Duflo $\tau_\nu$ \`a
l'alg\`ebre enveloppante de l'alg\`ebre de Lie $\g g^\C_\nu$,
d'espace vectoriel sous-jacent $\g g ^\C$ mais avec le crochet
d\'efini par $[x,y]_\nu=\nu[x,y]$. Pour un param\`etre r\'eel
$\hbar$ on notera $\g g_\hbar$ l'alg\`ebre de Lie r\'eelle
d'espace vectoriel sous-jacent $\g g$ mais avec le crochet
d\'efini par $[x,y]_\hbar=\hbar[x,y]$.
\alinea{II.2. Modules topologiquement libres convergents sur $(\Cal A,*)$}
Nous reprenons les notations du \S\ I.7. Rappelons que $\Cal A_0$
d\'esigne la sous-$\C[\nu]$-alg\`ebre de $\Cal A$ engendr\'ee par
$A$. Comme la s\'erie enti\`ere $a*b$ est un polyn\^ome en $\nu$
pour tout $a,b\in A$ on a~:
$$\Cal A_0=A[\nu].$$
Pour tout $\nu_0\in\C$ l'\'evaluation en $\nu_0$~:
$$\eqalign{\mop{ev}_{\nu_0}: \Cal A_0 &\longrightarrow (A,*_{\nu_0})\cr
        \sum_{k=0}^n\nu^ka_k          &\longmapsto   \sum_{k=0}^n\nu_0^ka_k \cr}$$
est un morphisme de $\C$-alg\`ebres.
\prop{II.2.1} Soit $R>0$, et soit pour tout $\nu_0\in D(0,R)$ une
repr\'esentation $\pi_{\nu_0}$ de l'alg\`ebre $(A,*_{\nu_0})$ dans
un espace vectoriel topologique localement convexe s\'epar\'e $M$.
On suppose que pour tout $a\in A$, $m\in M$ et $\nu_0\in D(0,R)$,
le vecteur $\pi_{\nu_0}(a)m$ est donn\'e par l'\'evaluation en
$\nu_0$ d'une s\'erie enti\`ere faiblement convergente de rayon
$\ge R$. Alors tout $\nu_0$ dans le disque de rayon $R$ induit une
repr\'esentation $\wt
{\pi_{\nu_0}}=\pi_{\nu_0}\circ\mop{ev}_{\nu_0}$ de $\Cal A_0$ dans
$M$, et $\Cal M=M[[\nu]]$ est alors un module topologiquement
libre faiblement convergent de rayon $\ge R$.
\smallskip
R\'eciproquement un module topologiquement libre faiblement
convergent de rayon $\ge R$ induit pour tout $\nu_0\in D(0,R)$ une
repr\'esentation $\pi_{\nu_0}$ de $(A,*_{\nu_0})$ dans $M$.
\dem
Soient $a,b\in A$ et $m\in M$. L'\'egalit\'e~:
$$\pi_{\nu_0}(a*_{\nu_0}b)m=\pi_{\nu_0}(a)\pi_{\nu_0}(b)m$$
pour tout $\nu_0\in D(0,R)$ implique l'\'egalit\'e entre s\'eries
formelles~:
$$\pi_\nu(a*b)m=\pi_\nu(a)\pi_\nu(b)m,$$
ce qui fait de $\Cal M=M[[\nu]]$ un module topologiquement libre.
Il est par construction faiblement convergent de rayon $\ge R$.
\smallskip
R\'eciproquement si $M$ est un espace vectoriel topologique
localement convexe s\'epar\'e et si $\Cal M=M[[\nu]]$ est un $\Cal
A$-module topologiquement libre faiblement convergent de rayon
$\ge R$, consid\'erons pour tout $\nu_0\in D(0,R)$ la
repr\'esentation $\wt{\pi_{\nu_0}}$ de $\Cal A_0$ dans $M$
donn\'ee par la proposition $I.7.1$. Si on note $\pi_{\nu_0}$ la
restriction de $\wt{\pi_{\nu_0}}$ \`a $A\subset \Cal A_0$, on a
imm\'ediatement pour tout $m\in M$~:
$$\eqalign{\pi_{\nu_0}(a)\pi_{\nu_0}(b)m &=\wt{\pi_{\nu_0}}(a)\wt{\pi_{\nu_0}}(b)m\cr
        &=\wt{\pi_{\nu_0}}(a*b)m    \cr
        &=\pi_{\nu_0}(a*_{\nu_0}b)m.\cr}$$
\qed

\alinea{II.3. Unitarit\'e}
On appelle repr\'esentation unitaire d'une alg\`ebre de Lie $\g g$
une repr\'esentation $\rho$ dans un espace pr\'ehilbertien telle
que les op\'erateurs $\rho(X)$ sont antihermitiens pour tout
$X\in \g g$. Le th\'eor\`eme suivant montre comment relier un
module topologiquement libre sur $\Cal A$ fortement unitaire et
convergent \`a une famille $\rho_\hbar$ de repr\'esentations
unitaires de l'alg\`ebre de Lie $\g g_\hbar$, lorsque $\hbar$
prend des valeurs r\'eelles.
\th{II.3.1}
Soit $R>0$, et soit $M$ un espace vectoriel topologique localement
convexe s\'epar\'e $M$ pr\'ehilbertien (i.e. muni d'un produit
scalaire hermitien). On suppose que $\Cal M=M[[\nu]]$ est un
module topologiquement libre faiblement convergent de rayon $\ge
R$. Soit pour tout $\nu_0\in D(0,R)$ la repr\'esentation
$\pi_{\nu_0}$ de l'alg\`ebre $(A,*_{\nu_0})$ dans $M$ associ\'ee
\`a $\Cal M$ par la proposition II.2.1. Alors~:
\smallskip 1°) l'\'egalit\'e~:
$$\rho_\hbar(X)=-i\pi_{i\hbar}(X)$$
d\'efinit une repr\'esentation $\rho_\hbar$ de l'alg\`ebre de Lie
$\g g_\hbar$ dans $M$ pour tout $\hbar\in D(0,R)$.
\smallskip
2°) Les deux propositions suivantes sont
\'equivalentes~:
\ssq
a) Pour tout $\hbar\in ]-R,R[$ la
repr\'esentation $\rho_\hbar$ est unitaire.
\ssq
b) Le module
topologiquement libre $\Cal M=M[[\nu]]$, muni du produit scalaire
\`a valeurs dans $\C[[\nu]]$ prolongeant naturellement celui de
$M$, est fortement unitaire.
\dem
On a pour $X,Y\in\g g$~:
$$\eqalign{[\rho_\hbar(X),\, \rho_\hbar(Y)] &=-[\pi_{i\hbar}(X),\,
\pi_{i\hbar}(Y)] \cr &=-i\hbar\pi_{i\hbar}([ X,Y])\cr
&=\hbar\rho_\hbar([X,Y]),\cr &=\rho_\hbar([X,Y]_\hbar),\cr}$$
d'o\`u la premi\`ere partie de la proposition. \ssq L'unitarit\'e
du module topologique $\Cal M=M[[\nu]]$ se traduit pour tout
$a\in\Cal A$ et $u,v\in\Cal M$ par l'\'egalit\'e entre s\'eries
formelles~:
$$<\pi_\nu(a)u,\, v>=<u,\,\pi_\nu(a^*)v>.$$
Comme nous avons suppos\'e que l'ind\'etermin\'ee $\nu$ est
imaginaire pure, cette \'egalit\'e se sp\'ecialise (au vu de la
proposition II.2.1) en tout param\`etre $\nu_0\in i]-R,R[$~: pour
tout $a\in A$ et $u,v\in M$ on a~:
$$<\pi_{\nu_0}(a)u,\, v>=<u,\,\pi_{\nu_0}(a^*)v>.$$
Comme $X^*=X$ pour tout $X\in \g g$ on voit imm\'ediatement que
les op\'erateurs $\pi_{i\hbar}(X)$ sont bien hermitiens pour
tout $\hbar\in ]-R,R[$. L'op\'erateur $\rho_\hbar(X )$ est donc
bien antihermitien, d'o\`u l'implication $b)\Rightarrow a)$. La
r\'eciproque est imm\'ediate.
\qed
\paragraphe{III. Exemples}
Nous explicitons ici la vari\'et\'e caract\'eristique et la
vari\'et\'e de Poisson caract\'eristique dans quelques exemples
Poisson-lin\'eaires. Le lemme suivant nous sera utile dans plusieurs de ces exemples ainsi qu'au chapitre IV~:
\lemme{III.0.1}
Soit $K$ un compact de $\R^n$ d'int\'erieur non vide, $m$ un entier non nul, $R>0$, et soit $(P_\nu)_{\nu\in D(0,R)}$ une famille d'op\'erateurs diff\'erentiels d'ordre $\le m$ dont les coefficients restreints \`a $K$ d\'ependent analytiquement de $\nu$. Alors pour tout $\varphi\in C^\infty_c(\R^n)$ \`a support inclus dans $K$ et pour toute distribution $T$ de support inclus dans $K$ la fonction  $\nu\mapsto <T,\, P_\nu(\varphi)>$ est analytique sur $D(0,R)$.
\dem
Soit $C^\infty_K(\R^n)$ l'espace des fonctions lisses sur $\R^n$ \`a support inclus dans $K$. On munit cet espace de la topologie de Fr\'echet d\'efinie par les seminormes~:
$$N_k(\varphi)=\mopl{sup}_{|\alpha|\le k}\mopl{sup}_{x\in K}|D^\alpha\varphi(x)|.$$
On consid\`ere le laplacien $\Delta$ sur $\R^n$. Il existe un entier $L$ tel que la distribution $(1-\Delta)^{-L}T$ soit une fonction continue sur $K$, car $T$ est d'ordre fini. On a alors par int\'egrations par parties~:
$$<T,P_\nu(\varphi)>=\int_K (1-\Delta)^{-L}T(x)(1-\Delta)^LP_\nu(\varphi)(x)\, dx.$$
On \'ecrit le d\'eveloppement en s\'erie enti\`ere~:
$$(1-\Delta)^LP_\nu(\varphi)(x)=P'_\nu(\varphi)(x)=\sum_{k\ge 0}\nu^kQ_k(\varphi)(x),$$
o\`u les $Q_k$ sont des op\'erateurs diff\'erentiels d'ordre $\le m+2L$. Les coefficients de $P'_\nu$ \'etant donn\'es sur $K$ par des s\'eries enti\`eres convergentes sur $D(0,R)$, il existe pour tout $r<R$ une constante $C$ telle que les coefficients de $Q_k$ sont tous major\'es sur $K$ en module par $Cr^{-k}$. On a donc~:
$$\mopl{sup}_{x\in K}|Q_k(\varphi)(x)|\le C'r^{-k}N_{m+2L}(\varphi),$$
d'o\`u la majoration~:
$$|\int_K(1-\Delta)^{-L}T(x)Q_k\varphi(x)\,dx|\le C'(\mop{Vol}K)N_{m+2L}(\varphi)\mopl{sup}_{x\in K}|(1-\Delta)^{-L}T(x)|r^{-k}.$$
Sur le disque de rayon $r$ l'expression $<T,\, P_\nu(\varphi)>$ est donc donn\'ee par la s\'erie enti\`ere convergente~:
$$<T,\, P_\nu(\varphi)>=\sum_{k\ge 0}\nu^k \int_K(1-\Delta)^{-L}T(x)Q_k\varphi(x)\,dx.$$
Ceci \'etant vrai pour tout $r<R$ on a bien convergence de cette s\'erie enti\`ere sur $D(0,R)$.
\qed
\alinea{III.1. Le groupe de Heisenberg}
Consid\'erons le groupe de Heisenberg $H_{2n+1}=~\R
^n\times\R ^n\times \R$ muni du produit
$$(x,y,z)\cdot (x',y',z')=(x+x',y+y', z+z'+{1\over 2}(xy'-x'y)).$$
Son alg\`ebre de Lie $\g h_{2n+1}$ est engendr\'ee par
$<X_1,...,X_n,Y_1,...,Y_n,Z>$ dont les crochets de Lie sont
donn\'es par~:
$$[X_i,X_j]=\delta_{ij}Z,~~i,j=1,...,n.$$

Pour $f_\l =\l Z^*+\displaystyle\sum_{i=1}^na_iX_i^*+
\displaystyle\sum_{i=1}^nb_iY_i^*\in {\g h^*_{2n+1}}$,
consid\'erons la repr\'esentation $\rho^\l$ associ\'ee \`a $f_\l$
par la m\'ethode des orbites.
\vskip 0,2cm
{\bf Cas 1~:} Si
$\l\not =0$ on peut prendre $f_\l =\l Z^*$, en effet, la dimension
de $\rho^\l$ est infinie et l'orbite associ\'ee \`a $f_\l$ sous
l'action coadjointe est $\O_{f_\l}=\{(\l,u,v),~ u,v\in~\R^n\}$.

Si on r\'ealise $\rho^\l$ \`a l'aide de la polarisation ${\g b}= <Y_1,Y_2,...,Y_n,Z>$,
alors $\rho^\l$ agit sur l'espace $L^2(\R^n)$. On notera aussi $\rho^\l$ sa diff\'erentielle,
qui se r\'ealise dans l'espace de Fr\'echet des vecteurs $C^\infty$, qui est ici l'espace de
Schwartz $\Cal S(\R^n)$ \cite{Ki}. Elle est donn\'ee par~:

$$\cases{\rho^\l (Z)=-i\l\cr\rho^\l(X_i)=-{\partial\over {\partial t_i}},~~i=1,...,n\cr
\rho^\l(Y_j)=i t_j,~~j=1,...,n.\cr}$$

Ainsi, la repr\'esentation $\rho_{\hbar}^\l$ est d\'efinie par les relations~:

$$\cases{\rho_{\hbar}^\l (Z)=-i\l\cr
\rho_{\hbar}^\l(X_i)=-{\partial\over {\partial
t_i}},~~i=1,...,n\cr \rho_{\hbar}^\l(Y_j)=i{\hbar}
t_j,~~j=1,...,n.\cr}$$
Le caract\`ere polynomial en $\hbar$ de ces expressions nous permet
d'utiliser la proposition II.2.1 et le th\'eor\`eme
II.3.1, ce qui fait de $\Cal M=\Cal S(\R^n)[[\nu]]$ un module topologiquement libre faiblement convergent (de rayon infini) fortement unitaire. La repr\'esentation de $\Cal {A}$
associ\'ee s'\'ecrit~:
$$\cases{\pi_\nu ^\l (Z)=\l\cr
\pi_\nu ^\l(X_i)=-i{\partial\over {\partial t_i}},~~i=1,...,n\cr
\pi_\nu ^\l(Y_j)=-i{\nu} t_j,~~j=1,...,n.\cr}$$
L'annulateur de $\pi^\lambda_\nu$ est donc engendr\'e par $Z-\lambda$, et donc~:
$$V\!A (\pi_\nu^\lambda)=\Omega_{f_\lambda}.$$
Par suite,
$$\cases{\pi_0^\l (Z)=\l\cr
\pi_0^\l(X_i)=-i{\partial\over {\partial t_i}},~~i=1,...,n\cr
\pi_0^\l(Y_j)=0,~~j=1,...,n.\cr}$$

On en d\'eduit alors que $\mop{Ann}(\pi_0^\l)$  est engendr\'e par
$<Z-\l,~Y_j, j=1,...,n>$.

Il vient alors que~:
$$\eqalign {V(\pi_\nu^\l)&=\{l\in{\g h_{2n+1}^*}/l(Z-\l)=0 ~
{\hbox{et}}~l(Y_j)=0 ,~j=1,...,n\}\cr &=\l Z^*\displaystyle
{\oplus_{i=1}^n} \R X^*_i=f_\l+{\g b}^{\bot}.\cr}$$

 Un calcul analogue  montre que si on r\'ealise $\rho^\l$ \`a l'aide de la polarisation

\noindent ${\g b'}=<X_1,...,X_n,Z>$, alors on a aussi que
$V(\pi_\nu^\l)=f_\l+{\g b'}^{\bot}$ et donc on voit clairement que
$V(\pi_\nu^\l)$ d\'epend de la r\'ealisation de $\pi^\l$.
\vskip
0,2cm {\bf Cas 2~:} Si $\l= 0$, alors, l'orbite $\O_{f_0}$ se
r\'eduit au point $\{f_0\}$ et la polarisation  associ\'ee \`a
$f_0$ est l'alg\`ebre de Lie $\g h_{2n+1}$ toute enti\`ere. Ainsi,
$$\rho^0(X)=-if_0(X),~~{\hbox{pour tout}}~ X\in \g h_{2n+1}.$$
Dans ce cas, on a que~: $\pi_{\nu}^0(X)=\pi_0^0(X),~{\hbox{pour
tout}}~ X\in \g h_{2n+1}$. On en d\'eduit alors que l'annulateur
de $\pi_\nu^0$ est l'id\'eal engendr\'e par
$\{X-f_0(X),~X\in \g h_{2n+1}\}$. Il s'ensuit alors que
$$\eqalign {V (\pi_\nu^0)&=\{l\in{\g h_{2n+1}^*}/l(X-f_0(X))=0,
~\forall X\in \g h_{2n+1}\}\cr &=\{f_0\},\cr}$$
et de m\^eme, $V\!A (\pi_\nu^0)=\{f_0\}$.
\goodbreak
\alinea{III.2. L'alg\`ebre de Lie filiforme de pas $n$}
Consid\'erons le groupe  de Lie nilpotent filiforme de pas $n$, $G_n$  d'alg\`ebre de Lie
$\g g_n$
de dimension $n+1$ muni d'une base de Jordan-H\"older  $(X_1,...,X_{n+1})$ avec
$$[X_{n+1},X_j]=X_{j-1},~~~j=2,...,n.$$
Soit $(X^*_1,\ldots ,X_{n+1}^*)$ la base duale de $\g g^*$. Le centre de $\g g$ est
engendr\'e par le vecteur $X_1$.
Soit $l=l_1X_1^*+\cdots +l_{n+1}X_{n+1}^*\in \g g^*$ avec $l_1\not =0$. Alors,
 $\g b (l)=$Vect$\{X_1,...,X_n\}$ est un id\'eal ab\'elien de $\g g_n$ qui polarise $l$.

\noindent Comme l'ensemble indice de Pukanszky est $\{2,n+1\}$,
sans perte de g\'en\'eralit\'e, on peut supposer que
$l_2=l_{n+1}=0$ (voir \cite{BBR}). La repr\'esentation unitaire et
irr\'eductible $\rho=\rho_l$ associ\'ee \`a $l$ se r\'ealise alors
 sur $L^2(\R )$. Sa diff\'erentielle est donn\'ee par~:
$$\cases{\displaystyle\rho (X_{n+1})=-{\partial\over{\partial t}}\cr
\rho(X_1)=-il_1\cr \rho(X_2)=itl_1\cr \displaystyle\rho
(X_3)=-i(l_3+{1\over 2}t^2l_1)\cr .\cr.\cr.\cr
\rho(X_n)=\displaystyle -i(l_n-l_{n-1}t+{1\over
2}l_{n-2}t^2+...+(-1)^{n-3}{{l_3}\over{(n-3)!}}
t^{n-3}+(-1)^{n-1}{{l_1}\over{(n-1)!}}t^{n-1}).\cr}$$ Il s'ensuit
alors que la repr\'esentation $\rho_\hbar$ est d\'etermin\'ee
par~:
$$\cases{\rho_\hbar (X_{n+1})= \displaystyle -{\partial\over{\partial t}}\cr
\rho_\hbar(X_1)=-il_1\cr
\rho_\hbar(X_2)=+i\hbar tl_1\cr
\rho_\hbar (X_3)=\displaystyle -i(l_3+{1\over 2}\hbar ^2t^2l_1)\cr
.\cr.\cr.\cr
\rho_\hbar(X_n)=\displaystyle -i(l_n-l_{n-1}\hbar t+{1\over 2}l_{n-2}\hbar^2 t^2+...
+{{l_3}\over{(n-3)!}}(-\hbar t) ^{n-3}+{{l_1}\over{(n-1)!}}(-\hbar t)^{n-1}).\cr}$$
Ces expressions \'etant polynomiales en $\hbar$ on peut encore appliquer la proposition II.2.1 et le th\'eor\`eme II.3.1, ce qui fait de $\Cal M=\Cal S(\R)[[\nu]]$ un module faiblement convergent fortement unitaire. On obtient~:
$$\cases{\pi_\nu (X_{n+1})= \displaystyle -i{\partial\over{\partial t}}\cr
\pi_\nu(X_1)=l_1\cr
\pi_\nu(X_2)=+i\nu tl_1\cr
\pi_\nu(X_3)=\displaystyle l_3-{1\over 2}\nu ^2t^2l_1\cr
.\cr.\cr.\cr
\pi_\nu(X_n)=\displaystyle l_n+i\nu l_{n-1}t-{1\over 2}\nu^2 l_{n-2} t^2+...
+(i\nu)^{n-3}{{l_3}\over{(n-3)!}}t^{n-3}+(i\nu)^{n-1}{{l_1}\over{(n-1)!}}t^{n-1}.\cr}$$
Faisant $\nu=0$ on obtient donc~:
$$\cases{\pi_0(X_{n+1})=\displaystyle -i{\partial\over{\partial t}}\cr
\pi_0(X_1)= l_1\cr
\pi_0(X_2)=0\cr
\pi_0 (X_3)= l_3\cr
.\cr.\cr.\cr
\pi_0 (X_n)= l_n.\cr}$$
En remarquant que Ann$(\pi_0)$ est engendr\'e par $<X_1- l_1,X_2,X_3- l_3,...,X_n- l_n>$,
il vient que
 pour $f=(f_1,...,f_{n+1})\in\g g^*$, $f\in V(\pi_\nu)$ si et seulement si
 $f\in l+\g b(l)^\bot$. D'autre part on voit que l'annulateur de $\pi_\nu$ est engendr\'e par
 les $v_k, k=1,\ldots ,n$, avec~:
$$v_k=X_k-l_k-{l_{k-1}\over l_1}X_2-{l_{k-2}\over 2l_1^2}X_2^2-\cdots
-{l_3\over (k-3)!l_1^{k-3}}X_2^{k-3}-{1\over (k-1)!l_1^{k-2}}X_2^{k-1}.$$
La vari\'et\'e de Poisson caract\'eristique $V\!A (\pi_\nu)$ est donc \'egale \`a l'orbite
coadjointe $\Omega_l$.
\alinea{III.3. L'alg\`ebre de Lie du groupe affine de la droite~:}
Consid\'erons le groupe de Lie compl\'etement r\'esoluble  $"aX+b"$ d\'efini par

\noindent $G=\{\pmatrix {a&b\cr 0&1\cr},~a,b\in~\R~~{\hbox
{et}}~a>0\}.$
 Son alg\`ebre de Lie est $\g g=\R X\oplus \R Y$,

\noindent o\`u $ X=\pmatrix {1&0\cr 0&0\cr} {\hbox{et}}~
Y=\pmatrix {0&1\cr 0&0\cr}$ et dont le crochet de Lie est donn\'e
par $[X,Y]=Y$.

\vskip 0,2cm
Ce groupe poss\`ede deux repr\'esentations unitaires et irr\'eductibles $\rho_+$ et $\rho_-$
associ\'ees respectivement aux formes lin\'eaires $Y^*$ et $-Y^*$.
Les diff\'erentielles de ces repr\'esentations sont d\'efinies par :

$$\cases{\rho_+(X)=-{d\over{dx}}\cr\rho_+(Y)=-ie^{-x}\cr}~~\hbox{et}~~
\cases{\rho_-(X)=-{d\over{dx}}\cr\rho_-(Y)=ie^{-x}.\cr}$$
Alors,
les expressions de $\rho_{+,{\hbar}}$ et
$\rho_{-,{\hbar}}$ sont donn\'ees respectivement par~:
$$\cases{\rho_{+,{\hbar}}(X)=-{d\over{dx}}\cr\rho_{+,{\hbar}}(Y)=-ie^{-{\hbar}x}\cr}
~~\hbox{et}~~
\cases{\rho_{-,{\hbar}}(X)=-{d\over{dx}}\cr\rho_{-,{\hbar}}(Y)=ie^{-{\hbar}x}.\cr}$$
On restreint ces diff\'erentielles \`a $M=C^\infty_K(\R)$ o\`u $K$ est un compact d'int\'erieur non vide. Appliquant d'abord le lemme III.0.1 puis la proposition II.2.1 et le th\'eor\`eme II.3.1, on fait de $\Cal M=M[[\nu]]$ un module faiblement convergent fortement unitaire (de rayon infini).
Les annulateurs des repr\'{e}sentations $\pi_{+\nu}$ et
$\pi_{-\nu}$ sont r\'eduits au singleton $\{ 0\}$ et par
cons\'equent, les ensembles $V\!A(\pi_{+,\nu})$ et $V\!A(\pi_{-,\nu})$
sont \'egaux \`a $\g g^*$. En prenant $\nu=0$, il vient que les
annulateurs de $\pi_{+,0}$ et $\pi_{-,0}$ sont respectivement
engendr\'es par $\langle Y-1\rangle$ et $\langle Y+1\rangle$. On
obtient alors que~:
$$\eqalign{V(\pi_{+,\nu})&=\{l=xX^*+yY^*\in\g g^*~/~\forall \f\in \mop{Ann}
\pi_{+,0},~\f (l)=0\}\cr
&=\{l=xX^*+yY^*\in\g g^*~/y-1=0\}=Y^*+\g b^\bot\cr}$$ o\`u $\g
b=\R Y$ est la polarisation asooci\'ee aux formes $Y^*$ et
$-Y^*$. De m\^eme, $$V(\pi_{-,\nu})=\{l=xX^*+yY^*\in\g
g^*~/y+1=0\}=-Y^*+\g b^\bot.$$
\vskip 0,2cm
\smallskip
{\sl Remarque}~: Pour ces deux repr\'esentations, la vari\'et\'e de Poisson caract\'eristique
co\"\i ncide avec l'adh\'erence de Zariski de l'orbite coadjointe associ\'ee.
\alinea{III.4. Un exemple r\'esoluble exponentiel de dimension 3:}
Soit $\g g$ l'alg\`ebre de Lie engendr\'ee par les
trois vecteurs $\{ A,X,Y\}$ dont les crochets de Lie sont donn\'es
par: $[A,X]=X-Y$, $[A,Y]=X+Y$ et soit $G=\exp \g g$.
 Alors, $G$ est un groupe  de Lie exponentiel non compl\`etement r\'esoluble.

Soit $f= xX^*+yY^*+aA^*\in \g g^*$. Si $x^2+y^2=0$, alors,
l'orbite de $f$ est r\'eduite au singleton $\{f\}$. Ainsi, le
calcul pr\'ec\'{e}dent de l'exemple III.1. montre que dans ce
cas, $V(\pi_{f,\nu})=\{ f\}$.

Dans le cas o\`u $x^2+y^2\not =0$, la sous-alg\`ebre $\g b$
engendr\'ee par $<X,Y>$ est une polarisation de $f$ v\'erifiant la
condition de Pukanszky. Soit alors $\chi_f$ le caract\`ere
d\'efini sur  $B=\exp \g b$ par $\chi_f(\exp U)=e^{-if(U)}$ et
$\rho_f= {\hbox{Ind}}_{B}^G \chi _f$. On sait \cite{ABLS} qu'il
existe un unique $\t\in [0,2\pi[$ tel que $\rho =\rho_\t
=\rho_{f_\t}$ o\`u $f_\t= \cos \t X^*+\sin \t Y^*$. L'orbite $\O$
associ\'ee \`a $\rho$ est param\'etris\'ee par
$$\O=\{sA^*+e^{-t}\cos (t+\t )X^*+e^{-t}\sin(t+\t )Y^*,~~s,t\in~\R~\}.$$

D'autre part, on a que~:

$$\cases{\rho(A)=-{d\over{dt}}\cr
\rho(X)=-ie^{-t}\cos (\t+t)\cr \rho(Y)=-ie^{-t}\sin (\t+t)\cr}
~~{\hbox{et}}~~ \cases{\rho_{\hbar}(A)=-{d\over{dt}}\cr
\rho_{\hbar}(X)=-ie^{-{\hbar}t}\cos (\t+{\hbar}t)\cr
\rho_{\hbar}(Y)=-ie^{-{\hbar}t}\sin (\t+{\hbar}t).\cr}$$
On se restreint comme dans l'exemple pr\'ec\'edent \`a $M=C^\infty_K(\R)$ o\`u $K$ est un compact d'int\'erieur non vide. Appliquant d'abord le lemme III.0.1 puis la proposition II.2.1 et le th\'eor\`eme II.3.1, on fait de $\Cal M=M[[\nu]]$ un module faiblement convergent fortement unitaire (de rayon infini). La repr\'esentation $\pi_\nu$ de $\Cal A$ s'\'ecrit~:
$$\cases{\pi_\nu(A)=-i{d\over{dt}}\cr
\pi_\nu(X)=e^{i\nu t}\cos (\t-i\nu t)\cr
\pi_\nu(Y)=e^{i\nu t}\sin (\t-i\nu t).\cr}
~~{\hbox{et}}~~
\cases{\pi_0(A)=-i{d\over{dt}}\cr
\pi_0(X)=\cos \t\cr \pi_0(Y)=\sin\t.\cr}$$
Il s'ensuit alors que l'annulateur de $\pi_\nu$ est r\'eduit \`a $\{0\}$, et donc~:
$$V\!A(\pi_\nu)=\g g^*.$$
La vari\'et\'e de Poisson caract\'eristique co\"\i ncide donc ici
aussi avec l'adh\'erence de Zariski de l'orbite coadjointe
associ\'ee. Pour le voir il suffit de se convaincre que la spirale
logarithmique est Zariski-dense dans le plan, en remarquant que
toute droite passant par l'origine intersecte cette spirale en une
infinit\'e de points. Par ailleurs l'annulateur de $\pi_0$ est
l'id\'eal engendr\'e par les deux g\'en\'erateurs $\{ X-\cos
\t,Y-\sin\t\}$. Soit $B=\exp \g b$. La repr\'esentation $\rho$
agit sur l'espace $L^2(G/B)$ qui est isomorphe \`a $L^2(\exp\R
A)$. On a alors,

$$\eqalign{V(\pi_\nu)&=\{ l\in \g g ^* ~/ l(X-\cos \t)=0 ~~{\hbox{et}}~~
l(Y-\sin \t)=0\} \cr
 &=f_\t + \R A^*=f_\t +\g b^\bot .\cr}$$
\goodbreak
\alinea{III.5. L'alg\`ebre de Lie du groupe diamant}
C'est l'alg\`ebre de Lie r\'eelle $\g g$ de dimension $4$ de base
$(H,P,Q ,E)$ avec les crochets~: $$[H,P]=-Q,\hskip 12mm
[H,Q]=P,\hskip 12mm [P,Q]=E,$$ les autres crochets \'etant nuls.
Notre r\'ef\'erence est M. Vergne dans \cite {BCD Chap.
VIII~\S~1.4}. Si $X=aH+bP+cQ+dE$ on calcule facilement la matrice
de $\mop{ad}X$ dans cette base~:
$$\mop{ad}X=\pmatrix{0 &0 &0 &0 \cr -c &0 &a &0 \cr b &-a &0 &0
\cr 0 &-c &b &0 \cr}.$$
On a \'egalement~:
$$(\mop{ad}X)^2=\pmatrix{0 &0 &0 &0 \cr ab &-a^2 &0 &0 \cr ac &0
&-a^2 &0 \cr b^2+c^2 &-ab &-ac &0 \cr},$$
ainsi que
l'\'egalit\'e~:
$$(\mop{ad}X)^3=-a^2(\mop{ad}X).$$
On en d\'eduit
l'expression explicite~:
$$\exp(\mop{ad}X)=I+{\sin a\over
a}\mop{ad}X+{1-\cos a\over a^2}(\mop{ad}X)^2.$$
Apr\`es
transposition et passage \`a l'inverse on en d\'eduit la matrice
de $\exp(\mop{ad}^*X)$ dans la base duale~:
$$\exp(\mop{ad}^*X)=\pmatrix{ 1 &c{\sin a\over a}+b{1-\cos a\over a}
&-b{\sin a\over a}+c{1-\cos a\over a} &(b^2+c^2){1-\cos a\over
a^2} \cr 0 &\cos a &\sin a &c{\sin a\over a}-b{1-\cos a\over a}
\cr 0 &-\sin a &\cos a &-b{\sin a\over a}-c{1-\cos a\over a} \cr 0
&0 &0 &1 \cr}.$$
L'action coadjointe de $\exp X$ sur un
\'el\'ement $\xi=\mu H^*+\beta P ^*+\gamma Q^*+\lambda E^*$
s'\'ecrit donc~:
$$\eqalign{\mop{Ad}^*(\exp X).\xi
&=\Big (\mu +(c{\sin a\over a}+b{1-\cos a\over a})\beta +(-b{\sin
a\over a}+ c{1-\cos a\over a})\gamma\cr &+(b^2+c^2){1-\cos a\over
a^2}\lambda \Big ) H^* \cr &+ \Big ((\cos a)\beta +(\sin a)\gamma
+(c{\sin a\over a}-b{1-\cos a\over a})\lambda\Big )P^*\cr &+ \Big
( (-\sin a)\beta +(\cos a)\gamma +(- b{\sin a\over a}-c{1-\cos
a\over a})\lambda\Big ) Q^*\cr &+\lambda E^*. \cr}$$ On note
$\xi_i$ la i-\`eme coordonn\'ee de $\mop{Ad}^*(\exp X).\xi$. La
derni\`ere coordonn\'ee $\xi_4$ est ici invariante sous l'action
coadjointe. Deux cas sont \`a consid\'erer~:
\smallskip {\sl
Premier cas}: $\xi_4=\lambda=0$. L'expression $\xi_2^2+\xi_3^2=
\beta^2+\gamma^2$ est invariante sur le sous-espace d\'efini par
$\xi_4= 0$. On voit alors que les orbites coadjointes
correspondantes sont les cylindres d'axe $H^*$ et chacun des
points de cet axe, suivant que $\xi_2^2 +\xi_3^2$ est strictement
positif ou s'annule. Ce sont exactement les orbites coadjointes du
quotient de $\g g$ par son centre (l'alg\`ebre de Lie du groupe
des d\'eplacements du plan)~: nous traitons cet exemple au
paragraphe suivant.
\medskip
{\sl Deuxi\`eme cas}~: $\xi_4=\lambda\not =0$. Faisant agir
$\exp Y$ sur $\xi$ avec $Y=b P+c Q$ on obtient~:
$$\eqalign{\mop{Ad}^*(\exp Y).\xi &= \bigl(\mu+c\beta+b\gamma+(b^
2+c^2)\lambda\bigr)H^* \cr &+(\beta +c\lambda)P^* \cr &+(\gamma-
b\lambda)Q^* \cr &+\lambda E^*. \cr}$$
En choisissant bien $b$ et $c$ on
peut donc se ramener au cas o\`u $\beta =\gamma=0$, ce que nous
supposerons. La formule explicite donnant $\mop{Ad}^*(\exp X).\xi$ se
simplifie alors~:
$$\eqalign{\mop{Ad}^*(\exp X).\xi
&=\Bigl(\mu+(b^2+c^2){1-\cos a\over a^2}\lambda\Bigr)H^* \cr &+
(c{\sin a\over a}-b{1-\cos a\over a})\lambda P^*\cr &+ (-b{\sin a\over
a}-c{1-\cos a\over a})\lambda Q^*\cr &+\lambda E^*. \cr}$$
On voit que l'on a~:
$$\eqalign{\xi_1-{{\xi_2}^2+{\xi_3}^2\over 2\xi_4} &=\mu \cr
\xi_4 &=\lambda.\cr}\eqno{(\III.5.1)}$$
Les orbites coadjointes
dans ce cas-l\`a sont donc les parabolo\"\i des de r\'evolution
$\Omega_{\lambda,\mu}$ d'axe $H^*$ donn\'es par les \'equations
III.5.1. \bigskip Soit $\varphi$ la fonction enti\`ere d'une
variable complexe d\'efinie par~: $$\varphi(z)={\mop{sh}z/2\over
z/2}.$$ Il est clair que la matrice de $\varphi^{1/2}(\mop{ad}X)$
est de la forme~:
$$\pmatrix{1 &0 &0 &0 \cr * &\varphi^{1/2}(ia) &0 &0 \cr * &*
&\varphi^{1/2}(ia) &0 \cr * &* &* &1 \cr},$$
d'o\`u avec les
notations du \S\ II.1~:
$$J(X)^{1/2}={\sin ({\nu a/2})\over \nu
a/2}.$$
L'isomorphisme de Duflo est donc donn\'e par~:
$$\tau={{\sin (\nu\partial_1 /2)}\over {\nu\partial_1 /2}}\circ\sigma,$$
o\`u $\partial_1$ d\'esigne l'op\'erateur de d\'eriv\'ee partielle
$\displaystyle{{\partial}\over {\partial\xi_1}}$ dans $\g g^*$, et
o\`u $\sigma$ d\'esigne la sym\'etrisation. Compte tenu des
\'equations III.5.1 on voit facilement que l'alg\`ebre $S(\g
g)^{\sg g}$ des polyn\^omes invariants sur $\g g^* $ est
engendr\'ee par $C_1$ et $C_2$, avec $C_1(\xi)=\xi_4$ et $C_2(\xi
)=\xi_2^2+\xi_3^2-2\xi_1\xi_4$. L'orbite coadjointe
$\Omega_{\lambda,\mu}$ d\'efinie par III.5.1 est donn\'ee de fa\c
con \'equivalente par les \'equations~:
$$C_1(\xi)=\lambda,\hskip 12mm C_2(\xi)=-2\lambda\mu.
\eqno{(\III.5.2 )}$$
L'op\'erateur $J(D)^{1/2}$ agit par l'identit\'e sur
ces deux g\'en\'erateurs, et donc l'isomorphisme de Duflo $\tau$ se
ram\`ene \`a appliquer la sym\'etrisation sur les deux g\'en\'erateurs.
Nous noterons encore un peu abusivement $C_1$ et $C_2$ les deux
Casimirs $\tau (C_1)$ et $\tau (C_2 )$. On a explicitement~:
$$C_1=E,\hskip 12mm C_2=P^2+Q^2-2EH.$$
Nous utilisons maintenant les notations du \S\ II. Pour tout
r\'eel $\hbar$ il existe une famille de repr\'esentations
unitaires irr\'eductibles $ \rho_{\lambda,\mu;\hbar}$ de
$G_\hbar=\exp \g g_\hbar$, dont les diff\'erentielles sont
d\'efinies sur l'espace de Schwartz $\Cal S(\R)$ (muni du produit
scalaire de $L^2(\R)$) par~:
$$\eqalign{\rho_{\lambda,\mu;\hbar}(H) &=-i(-{1\over
2\lambda}{d^2\over dx^2}+{1\over 2}\lambda\hbar^2x^2+\mu) \cr
\rho_{\lambda,\mu;\hbar}(P) &=-{d\over dx} \cr
\rho_{\lambda,\mu;\hbar}(Q) &=i\lambda\hbar x \cr
\rho_{\lambda,\mu;\hbar}(E) &=-i\lambda. \cr}$$
Ces
repr\'esentations sont unitairement \'equivalentes \`a celles que l'on
obtient par induction holomorphe \`a partir du point $\mu
H^*+\lambda E^* $ et de la polarisation complexe $\g h$
engendr\'ee par $H,E,P+iQ$ (\cite {BCD} \S\ V.4 et VIII.1.4.4).
Les expressions sont ici polynomiales en $\hbar$. Posant $\nu=i\hbar$ et appliquant la proposition II.2.1 et le th\'eor\`eme II.3.1 on en d\'eduit
une repr\'esentation unitaire $\pi_{\lambda,\mu; \nu}$ de
l'alg\`ebre d\'eform\'ee $\Cal A$ dans le module topologiquement
libre $\Cal S(\R)[[\nu]]$, faiblement convergent fortement unitaire~:
$$\eqalign{\pi_{\lambda,\mu;\nu}(H)
&=-{1\over 2\lambda}{d^2\over dx^2 }-{1\over 2}\lambda\nu^2x^2+\mu
\cr \pi_{\lambda,\mu;\nu}(P) &=-i{d\over dx} \cr
\pi_{\lambda,\mu;\nu}(Q) &=i\lambda\nu x \cr
\pi_{\lambda,\mu;\nu}(E) &=\lambda. \cr}$$
En annulant le
param\`etre $\nu$ on voit que l'annulateur $\mop{Ann}\pi_{
\lambda,\mu;0}$ est l'id\'eal de $S(\g g)$ engendr\'e par
$E-\lambda$, $Q $ et $C_2+2\lambda\mu$. On en d\'eduit que la
vari\'et\'e caract\'eristique $V(\pi_{\lambda,\mu;\nu})$ est
d\'efinie par les \'equations~:
$$\xi_2^2+\xi_3^2-2\xi_1\xi_4=-2\lambda\mu,\hskip 12mm
\xi_4=\lambda, \hskip 12mm \xi_3=0.$$ La vari\'et\'e
caract\'eristique $V(\pi_{\lambda,\mu;\nu})$ est donc une
g\'en\'eratrice du parabolo\"\i de de r\'evolution
$\Omega_{\lambda,\mu}$. La non-existence de polarisations
r\'eelles se traduit ici par le fait que la vari\'et\'e
caract\'eristique n'est pas un sous-espace affine (cf. \S\ IV.2).
\medskip
L'annulateur
$\mop{Ann}\pi_{\lambda,\mu;\nu}$ est quant \`a lui engendr\' e
dans $\Cal A$ par $C_1-\lambda$ et $C_2+2\lambda\mu$. L'id\'eal
$\mop{Ann}\pi_{\lambda,\mu;\nu}/\nu
\mop{Ann}\pi_{\lambda,\mu;\nu}$ est donc engendr\'e dans $S(\g g)$
par $C_1-\lambda$ et $C_2+2\lambda\mu$ \'egalement. On en d\'eduit
que la vari\'et\'e de Poisson caract\'eristique $V\!A(
\pi_{\lambda,\mu;\nu})$ co\"\i ncide avec l'orbite coadjointe
$\Omega_{\lambda,\mu}$.
\alinea{III.6. L'alg\`ebre de Lie du groupe des d\'eplacements du plan}
C'est l'alg\`ebre de Lie r\'eelle $\g g$ de dimension $3$ de base
$(H,P,Q )$ avec les crochets~: $$[H,P]=-Q,\hskip 12mm [H,Q]=P,$$
les autres crochets \'etant nuls. Nous nous r\'ef\'erons encore
\`a l'article de M. Vergne \cite {BCD Chap. VIII~\S~1.3}, ainsi
qu'\`a \cite {AM}. Si $X=aH+bP+cQ$ on calcule facilement la
matrice de $\mop{ad}X$ dans cette base~:
$$\mop{ad}X=\pmatrix{ 0
&0 &0 \cr -c &0 &a \cr b &-a &0 \cr}.$$
On a \'egalement~:
$$(\mop{ad}X)^2=\pmatrix{0 &0 &0 \cr ab &- a^2 &0 \cr ac &0 &-a^2
\cr },$$
ainsi que l'\'egalit\'e~:
$$(\mop{ad}X)^3=-
a^2(\mop{ad}X).$$
On en d\'eduit la m\^eme expression explicite
que dans l'exemple pr\'ec\' edent~:
$$\exp(\mop{ad}X)=I+{\sin
a\over a}\mop{ad}X+{1-\cos a\over a^2}(\mop{ad}X)^2.$$
Apr\`es
transposition et passage \`a l'inverse on en d\'eduit la matrice
de $\exp(\mop{ad}^*X)$ dans la base duale~:
$$\exp(\mop{ad}^*X)=\pmatrix{ 1 &c{\sin a\over a}+b{1-\cos a\over a}
&-b{\sin a\over a}+c{1-\cos a\over a} \cr 0 &\cos a &\sin a \cr 0 &-
\sin a &\cos a }.$$
L'action coadjointe de $\exp X$ sur un \'el\'ement
$\xi=\mu H^*+\beta P ^*+\gamma Q^*$ s'\'ecrit donc~:
$$\eqalign{\mop{Ad}^*(\exp X).\xi &=\Bigl(\mu +(c{\sin a\over a}+ b{1-
\cos a\over a})\beta +(-b{\sin a\over a}+c{1-\cos a\over a})\gamma\Bigr)H^* \cr &+
 \Bigl( (\cos a)\beta +(\sin a)\gamma\Bigr)P^*\cr &+ \Bigl(
(-\sin a)\beta +(\cos a)\gamma\Bigr)Q^*.\cr }$$
Les orbites coadjointes
sont donc les points $\mu H^*$ et les cylindres $ \Omega_r,r>0$ d'axe
$H^*$ et de rayon $r$, d\'efinis par l'\'equation $\xi_2^2+\xi_3^2=r^2$.
\ssq
L'isomorphisme de Duflo est encore donn\'e par~:
$$\tau={\sin ({\nu \p_1/2})\over \nu \p_1 /2}\circ\sigma,$$
o\`u $\p_1$
d\'esigne l'op\'erateur de d\'eriv\'ee partielle
$\displaystyle{{\p}\over {\p \xi_1}}$
 dans $\g g^*$, et o\`u $\sigma$ d\'esigne
la sym\'etrisation. L'alg\`ebre des invariants admet
$C~:\xi\mapsto\xi_2^2+ \xi_3^2$ comme seul g\'en\'erateur, et
l'isomorphisme de Duflo consiste encore \`a appliquer l'identit\'e
sur ce g\'en\'erateur. Le Casimir $\tau( C)$ que nous notons
encore $C$ s'\'ecrit~:
$$C=P^2+Q^2.$$ On consid\`ere pour $r>0,
\hbar\not =0$ et $0\le \lambda<1$ la famille
$\rho_{r,\lambda;\hbar}$ de repr\'esentations unitaires
irr\'eductibles de $G_\hbar=\exp \g g_\hbar$ suivantes~: la repr\'esentation $\rho_{r,\lambda;\hbar}$ agit sur l'espace $\Cal H_\lambda$ des classes de fonctions $2\pi$-pseudo-p\'eriodiques $\varphi$, telles que~:
$$\varphi(t+2\pi)=e^{2i\pi\lambda}\varphi(t),\hskip 12mm
\int_0^{2\pi}|\varphi(t)|^2\,dt<+\infty.$$
Sa diff\'erentielle est d\'efinie par~:
$$\eqalign{\rho_{r,\lambda;\hbar}(H) &=-\hbar{d\over d\theta} \cr
\rho_{r,\lambda;\hbar}(P) &=ir\sin (\theta) \cr
\rho_{r,\lambda;\hbar}(Q) &=ir\cos (\theta) \cr}$$
sur l'espace $\Cal H_\lambda^\infty$
des vecteurs $C^\infty$, qui est ici l'espace des fonctions $C^\infty$ de $\Cal H_\lambda$. Ces expressions sont polynomiales en $\hbar$. Posant $\nu=i\hbar$
et appliquant la proposition II.2.1 et le th\'eor\`eme II.3.1 on a donc une famille de
repr\'esentations unitaires de $\Cal A$ dans le module topologiquement libre faiblement convergent fortement unitaire $\Cal H_\lambda^\infty[[\nu]]$~:
$$\eqalign{\pi_{r,\lambda;\nu}(H) &=-\nu{d\over d\theta} \cr
\pi_{r,\lambda;\nu}(P) &=-r\sin (\theta) \cr
\pi_{r,\lambda;\nu}(Q) &=-r\cos (\theta). \cr}$$
En annulant $\nu$ on voit que l'annulateur
$\mop{Ann}\pi_{r,\lambda;0}$ est engendr\'e par $H$ et par $P^2+Q^2-
r^2$. La vari\'et\'e caract\'eristique $V(\pi_{r,\lambda;\nu})$ est donc
le cercle de cote nulle dans $\Omega_r$ d\'efini par les \'equations
$\xi_2^2+\xi_3^2=r^2$ et $\xi_1=0$. Elle est donc ind\'ependante du
param\`etre suppl\'ementaire $\lambda$. L'annulateur
$\mop{Ann}\pi_{r,\lambda;\nu}$ est quant \`a lui engendr\'e dans $\Cal
A$ par $C-r^2$. L'id\'eal $\mop{Ann}\pi_{r,\lambda;\nu}/\nu \mop{Ann}\pi_{r,\lambda;\nu}$ est
 donc \'egalement engendr\'e par $C-r^2$ dans $S(\g g)$.
La vari\'et\'e de Poisson caract\'eristique
$V\!A(\pi_{r,\lambda;\nu})$ est donc \'egale \`a l'orbite
coadjointe $\Omega_r$.
\ssq
On consid\`ere maintenant l'espace $\Cal H_\lambda^\omega$ des vecteurs analytiques. Ce sont
les fonctions $\varphi$ enti\`eres sur $\R$ telles que
$\varphi(t+2\pi)=e^{2i\pi\lambda}\varphi(t)$.
Soit $R$ l'op\'erateur sur
$C^\omega(\R)[[\nu]]$ d\'efini par~:
$$R\varphi(\theta)=\nu\varphi(\nu\theta).$$
Cet op\'erateur est
injectif. Soit $\widetilde C_\lambda=R(\Cal H_\lambda^\omega[[\nu]])$, et soit
$\widetilde\pi_{r,\lambda;\nu}$ la repr\'esentation de $\Cal A$
sur $\widetilde C_\lambda$ d\'efinie par~:
$$\widetilde\pi_{r,\lambda;\nu}(X)=R\circ\pi_{r,\lambda;\nu}\circ
R\inver.$$
Les deux repr\'esentations sont \'equivalentes par
construction, et on a~:
$$\eqalign{\widetilde\pi_{r,\lambda;\nu}(H)
&=-{d\over d\theta} \cr \widetilde\pi_{r,\lambda;\nu}(P) &=-r\sin
(\nu\theta) \cr \widetilde\pi_{r,\lambda;\nu}(Q) &=-r\cos (\nu\theta).
\cr}$$
L'annulateur de $\widetilde\pi_{r,\lambda;0}$ est engendr\'e par
$P$ et par $Q+r$. La vari\'et\'e caract\'eristique
$V(\widetilde\pi_{r,\lambda;\nu})$ est donc cette fois-ci la
g\'en\'eratrice du cylindre $\Omega_r$ d\' efinie par les \'equations
$\xi_3=-r$ et $\xi_2=0$.
\alinea{III.7. Le cas semi-simple~: modules de Verma}
Soit $\g g_1$ une alg\`ebre de Lie semi-simple complexe. La forme de
Killing d\'efinie par~:
$$(X,Y)=\mop{Tr}(\mop{ad}X\circ\mop{ad}Y)$$
est non d\'eg\'en\'er\'ee et invariante par l'action adjointe. Elle induit
donc un isomorphisme lin\'eaire $\kappa$ de $\g g_1$ sur son dual $\g
g _1^*$, d\'efini par~:
$$\kappa(X)=(X,-),$$
qui entrelace les repr\'esentations adjointe et coadjointe. Soit
$H$ un \'el\'ement semi-simple de $\g g_1$, soit $\g h_1$ une
sous-alg\`ebre de Cartan contenant $H$, soit $\Delta$ le syst\`eme
de racines associ\'e, et soit $\Delta_+$ l'ensemble des racines
positives provenant du choix d'un ordre sur $\g h_1^*$. On
d\'esignera par $W$ le groupe de Weyl associ\'e \`a ces donn\'ees.
\ssq La forme de Killing restreinte \`a $\g h_1$ est non
d\'eg\'en\'er\'ee. Soit $\lambda=\kappa\inver(H)\in\g g_1^*$. La
d\'ecomposition radicielle de $\g g_1$ s'\'ecrit~:
$$\g g_1=\g n_{1-}\oplus\g h_1\oplus\g n_{1+},$$
o\`u $\g
n_{1\pm}=\bigoplus_{\alpha\in\pm\Delta_+}\g g_1^\alpha.$ Soit
$\delta\in\g h_1^*$ la demi-somme des racines positives. La sous-alg\`ebre
de Cartan $\g h_1$ est orthogonale \`a $\g n_{1+}$ et $\g
n_{1-}$ pour la forme de Killing, ce qui permet de montrer que
$\kappa\inver(H)\in\g g _1^*$ s'annule sur $\g n_{1\pm}$. On
identifiera donc $\g h_1^*$ \`a l'orthogonal de $\g n_{1-}\oplus\g
n_{1+}$ dans $\g g_1^*$ en prolongeant trivialement \`a $\g n_{1-
}\oplus\g n_{1+}$ les formes lin\'eaires sur $\g h _1$. Si $H$ est
r\'egulier le stabilisateur de $\lambda$ pour l'action coadjointe est
\'egal \`a $\g h_1$, et la sous-alg\`ebre r\'esoluble $\g b_ 1=\g
h_1\oplus\g n_{1+}$ est une polarisation r\'esoluble en $\lambda$.
\ssq
Le cadre ci-dessus s'applique aux alg\`ebres de Lie $(\g
g_\nu)_{\nu\in\C -\{0\}}$ (avec les notations du \S ~II.2), qui ont toutes
le m\^eme espace vectoriel sous-jacent que nous noterons $\g g$. Nous
continuerons \`a identifier $\g g$ et son dual avec la forme de Killing de
$\g g=\g g_1$ (ind\'ependamment de $\nu$) et non pas avec celle de
$\g g_\nu$. La d\'ecomposition radicielle~:
$$\g g_\nu=\g n_{\nu-}\oplus\g h_\nu\oplus\g n_{\nu+}$$
est ind\'ependante de $\nu$ en ce
qui concerne les espaces vectoriels sous-jacents que nous noterons $\g
n_-$, $\g h$ et $\g n_+$. De m\^eme, on note $\g b_+=\g h\oplus \g n_+$ et~:
$$\g g^\alpha=\{X\in\g g,\, \forall
A\in\g h, [A,X]_\nu=\nu\alpha(A)X \}.$$
Le syst\`eme de racines
associ\'e \`a $\g g_\nu$ est donc $\nu\Delta$, et $\g g^\alpha$ co\"\i
ncide pour tout $\nu$ avec le sous-espace radiciel de $\g g_\nu$
correspondant \`a $\nu\alpha$. Soit $M_\lambda^\nu$ le module de
Verma associ\'e \`a $\lambda$ pour $\g g_\nu$. On a~:
$$M_\lambda^\nu=\Cal U(\g g_\nu)\mathop{\otimes}\limits_{\Cal U(\sg
b_{\nu +})}\C,$$ o\`u un \'el\'ement $A+Y$ de $\g b_{\nu +}=\g
h_\nu\oplus\g n_{\nu+}$ agit sur $\C$ via la multiplication par
$(\lambda- \nu\delta)(A)$, et o\`u l'action de $\Cal U(\g g_\nu)$
est donn\'ee par la multiplication \`a gauche. Lorsqu'on regarde
$\nu$ comme une ind\'etermin\'ee ce module est topologiquement
libre sur $\Cal A=\Cal U_\nu(\g g)$~: il s'identifie en effet \`a $\Cal U_\nu(\g n_-)$, et
donc \`a $S(\g n _-)[[\nu]]$ via sym\'etrisation. On a la
d\'ecomposition en sous-espaces de poids sous l'action de $\g
h_\nu$~:
$$M_\lambda^\nu=\bigoplus_{\beta\in
Q_+}(M_\lambda^\nu)_{\lambda-\nu\delta-\nu\beta},$$
o\`u $Q_+$
d\'esigne l'ensemble des combinaisons lin\'eaires \`a coefficients entiers
positifs d'\'el\'ements de $\Delta_+$. En annulant $\nu$ on voit que tout
$A\in \g h_0$ agit sur $M_\lambda^0$ par multiplication par
$\lambda(A)$, que l'action de $\g n_{0+}$ est triviale sur
$M_\lambda^0$ , et que l'action de $\g n_{0-}$ est fid\`ele. L'annulateur
de $M_\lambda ^0$ dans $S(\g g)$ est donc l'id\'eal engendr\'e par $\g
n_{+}$ et par les $A-\lambda(A),\, A\in\g h$. On en d\'eduit~:
$$V(M_\lambda^\nu)=\lambda+\g b_+^\perp.$$
Soit
$\chi_\lambda^\nu~:Z(\Cal U(\g g_\nu))\to \C$ le caract\`ere
central de $M_\lambda^\nu$. Soit $\gamma_\nu ~: Z(\Cal U(\g
g_\nu))\widetilde{\longrightarrow} S(\g h)^W$ l'isomorphisme
d'Harish-Chandra de $\g g_\nu$ \cite{Dix~\S~7.4}. En consid\'erant
$S(\g h)^W$ comme l'ensem\-ble des polyn\^omes $W$-invariants sur
$\g h^*$ on a alors \cite{Dix~\S~7.4.6}~:
$$\chi_\lambda^\nu(u)=\gamma_\nu(u)(\lambda).$$
L'isomorphisme de
Duflo $\tau_\nu~:S(\g g)^{\sg g_\nu}\to Z(\Cal U(\g g_\nu))$ s'obtient
en composant $\gamma_\nu\inver$ \`a gauche par la restriction \`a
$\g h^*$ (isomorphisme de Chevalley). Pour tout $v\in S(\g g)^{\sg
g_\nu} $ on a donc finalement~:
$$\chi_\lambda^\nu(\tau_\nu(v))=v(\lambda).$$
\qquad L'annulateur du module de Verma $M_\lambda^\nu$ est l'id\'eal
bilat\`ere de $\Cal U(\g g_\nu)$ engendr\'e par
$\mop{Ker}\chi_\lambda^\nu$ (\cite{Dix}, Th\'eor\`eme~8.4.3). Via
l'isomorphisme de Duflo c'est donc l'id\'eal bilat\`ere de
$(S(\g g),*_\nu)$ engendr\'e par $\{v-v(\lambda),\,v\in S( \g g)^{\sg g_\nu}\}$.
\ssq
Consid\'erons maintenant $\nu$ comme une ind\'etermin\'ee et $M_\lambda^\nu$ comme module topologiquement libre sur $\Cal A$. L'id\'eal $\mop{Ann}M_\lambda^\nu /
(\mop{Ann}M_\lambda^\nu\cap \nu\Cal A)$ de $S (\g g)$ est donc l'id\'eal engendr\'e
par $\{v-v(\lambda),\,v\in S(\g g)^ {\sg g_\nu}\}$. On en d\'eduit la
vari\'et\'e de Poisson caract\'eristique~:
$$V\!A(M_\lambda^\nu)=\{\xi\in\g g^*,\, v(\xi)=v(\lambda)\hbox{ pour
tout }v\in S(\g g)^{\sg g_\nu} \}.$$
Lorsque $\lambda$ (c'est-\`a-dire $H$)
est r\'egulier, c'est l'orbite coadjointe de $\lambda$. Lorsque
$\lambda=0$ c'est le c\^one nilpotent.
\alinea{III.8. Modules de Verma et r\'ealit\'e}
Le module de Verma $M_\lambda^\nu$ admet une forme bilin\'eaire
sym\'etrique naturelle ayant de bonnes propri\'et\'es de
covariance~: la forme de Shapovalov \cite {Sha}, \cite{D2}. En la
modifiant pour la rendre sesquilin\'eaire hermitienne nous
pourrons appliquer les r\'esultats du \S\ I.6. \ssq On garde les
notations du \S\ III.7. Soit
$(X_{-\alpha},H_{\alpha},X_{\alpha})_{\alpha\in\Delta_+}$ une base
de Chevalley de $\g g_1$. Alors
$(-iX_{-\alpha},-iH_{\alpha},-iX_{\alpha})_{\alpha\in\Delta_+}$
est une base de Chevalley de $\g g_i$ (o\`u $i=\sqrt{-1}$). Soit
$\g g_{i,\R}$ la forme r\'eelle d\'eploy\'ee de $\g g_i$
associ\'ee. Elle est d\'efinie comme l'espace vectoriel r\'eel
engendr\'e par les $-iX_{-\alpha},-iH_\alpha, -iX_\alpha$. Ceci
d\'efinit une forme r\'eelle $\g g_\R$ de l'espace vectoriel
sous-jacent $\g g$, qui est \`a son tour automatiquement une forme
r\'eelle d\'eploy\'ee $\g g_{\nu,\R}$ de l'alg\`ebre de Lie $\g
g_\nu$ pour tout $\nu$ imaginaire pur. La conjugaison respecte les
sous-espaces radiciels, donc a fortiori les trois composantes de
la d\'ecomposition~:
$$\g g_\nu=\g n_{\nu-}\oplus\g h_\nu\oplus\g n_{\nu+}.$$
La conjugaison sur $\g g$ (resp. sur $\g h$) par rapport \`a cette
forme r\'eelle induit une conjugaison sur le dual $\g g^*$ (resp.
$\g h^*$) gr\^ ace \`a la formule~:
$$\overline \xi (X):=\overline{\xi(\overline X)}.$$
Nous \'etendons la conjugaison par multiplicativit\'e \`a
l'alg\`ebre enveloppante $\Cal U(\g g_\nu)$. Cette alg\`ebre
enveloppante admet (gr\^ace au th\'eor\`eme de
Poincar\'e-Birkhoff-Witt) la d\'ecomposition~:
$$\Cal U(\g g_\nu)=\Cal U(\g h_\nu)\oplus\bigl(\g n_{\nu-}\Cal U(\g g_\nu)
    +\Cal U(\g g_\nu)\g n_{\nu+}\bigr).$$
Soit $P_\nu:\Cal U(\g g_\nu)\to\Cal U(\g h_\nu)$ la projection
correspondant \`a cette d\'ecomposition. D'apr\`es ce qui
pr\'ec\`ede cette projection est compatible avec la conjugaison,
c'est-\`a-dire~:
$$P(\overline u)=\overline{P(u)}\hskip 6 mm \hbox{pour tout }u\in \Cal U(\g g_\nu).$$
Nous d\'efinissons la {\sl transposition} sur $\Cal U(\g g_\nu)$ de la fa\c con suivante~:
$$^t\!X_\alpha=X_{-\alpha},\hskip 10mm
^t\!X_{-\alpha}=X_{\alpha},\hskip 10mm
^t\!H_{\alpha}=H_{\alpha},$$ et on \'etend cette transposition en
un anti-automorphisme de $\Cal U(\g g_\nu)$. Il est imm\'ediat de
voir que la transposition commute avec la conjugaison. Nous
consid\'ererons sur $\Cal U(\g g_\nu)$ l'involution
(semi-lin\'eaire) d\'efinie par~:
$$a^*=^t\!\overline a.$$
Cette involution (restreinte \`a $\g g$) d\'efinit une nouvelle
conjugaison sur $\g g$, et donc une nouvelle forme r\'eelle $\g
g^\R$, qui est une forme r\'eelle {\sl compacte\/} de chacune des
alg\`ebres de Lie $\g g_\hbar$ o\`u $\hbar=-i\nu$ est un r\'eel
non nul (cf. \cite {H} \S\ III.6). Remarquons que l'on a~:
$$\g h\cap\g g^\R=i\g h\cap\g g_{\R}.$$
Nous noterons $\g h^{\R}$ cette intersection et nous noterons
$\tau$ la conjugaison (dans $\g g$ ou dans $\g h$ ou dans leurs
duaux respectifs) par rapport \`a cette nouvelle forme r\'eelle.
Supposons alors que $\lambda\in\g h^*$ est r\'eel, c'est-\`a-dire
$\tau(\lambda) =\lambda$. L'involution $*$ sur $\Cal U(\g g_\nu)$
co\"\i ncide, via l'isomorphisme de Duflo, avec l'involution
donn\'ee au \S\ I, la conjugaison \'etant $\tau$.
\ssq
Plut\^ot
que de consid\'erer le module de Verma $M_\nu^\lambda$ comme dans
le \S\ III.7 nous consid\'ererons le module de Verma
$M_\nu^{\lambda+\nu\delta}$. Consid\'erons maintenant $\nu$ comme une ind\'etermin\'ee~: en tant que modules sur l'alg\`ebre
d\'eform\'ee $\Cal A$ ils co\"\i ncident \`a $O(\nu)$ pr\`es donc leurs
vari\'et\'es caract\'eristiques et leurs vari\'et\'es de Poisson
caract\'eristiques sont les m\^emes. Soit
$e_{\lambda+\nu\delta}^\nu$ le vecteur de plus haut poids de
$M_{\lambda+\nu\delta}^\nu$. \ssq Supposons que $\lambda$ soit
r\'eel. Nous d\'efinissons alors une forme sesquilin\'eaire sur le
module de Verma $M_{\lambda+\nu\delta}^\nu$ par la formule~:
$$<m,n>_\nu=<a.e_{\lambda+\nu\delta}^\nu,b.e_{\lambda+\nu\delta}^\nu>_\nu=P_\nu(a^*b)(\lambda).$$
C'est une version hermitienne de la forme de Shapovalov de $M_{\lambda+\nu\delta}^\nu$.
\prop{III.8.1}
La repr\'esentation $\pi_\nu$ de $\Cal A=\Cal U(\g g_\nu)$ dans
$M_{\lambda+\nu\delta}^\nu$ v\'erifie pour tout $m,n\in
M_\lambda^\nu$ et pour tout $u\in\Cal U(\g g_\nu)$~:
$$<um,n>_\nu=<m,u^*n>_\nu.$$
\dem Soient $a,b\in\Cal U(\g g_\nu)$ tels que $m=a.e_\lambda^\nu$
et $n=b.e_\lambda^\nu$. La proposition d\'ecoule imm\'ediatement
de l'\'egalit\'e~:
$$<ua.e_{\lambda+\nu\delta}^\nu,\,b.e_{\lambda+\nu\delta}^\nu>=
<a.e_{\lambda+\nu\delta}^\nu,\,u^*b.e_{\lambda+\nu\delta}^\nu>.$$
\qed La repr\'esentation associ\'ee est donc une
$*$-repr\'esentation. De plus la forme quotient d\'efinie par
$<-,->_\nu$ en $\nu=0$ est hermitienne non d\'eg\'en\'er\'ee.
Autrement dit le module de Verma $M_{\lambda+\nu\delta}^\nu$ est
fortement pseudo-unitaire. Nous pouvons donc appliquer les
r\'esultats du \S\ I.6. Les vari\'et\'es caract\'eristiques
$V(\pi_\nu)$ et $V\!A(\pi_\nu)$ sont donc d\'efinies sur le corps
des r\'eels. Elles sont donn\'ees par les m\^emes formules
explicites que dans le \S\ III.7.
\smallskip
{\sl Remarque}~: Soit $\lambda\in\g h^*$. Alors, en utilisant le
th\'eor\`eme 7.6.24 de \cite {Dix}, il est facile de voir que le
module de Verma $M_{\lambda+\nu\delta}^\nu$ est simple (et donc
que la forme de Shapovalov est non-d\'eg\'en\'er\'ee) sauf
\'eventuellement sur un ensemble discret d\'enombrable de valeurs
de $\nu$. Par ailleurs, sur la question de l'unitarisabilit\'e de
certains modules de plus haut poids, voir \cite {EHW} et \cite
{Jak}.
\paragraphe{IV. Le cas r\'esoluble exponentiel}
\qquad Dans le cas r\'esoluble exponentiel, nous explicitons dans ce
paragraphe, en utilisant la construction de \cite {Pe2}, la
vari\'et\'e caract\'eristique d'une repr\'esentation $\pi_\nu$
associ\'ee de mani\`ere naturelle \`a une repr\'esentation
monomiale unitaire quelconque induite \`a partir d'une polarisation
r\'eelle.
\alinea{IV.1. Rappels sur la m\'ethode des orbites pour les groupes exponentiels}
\qquad L'adaptation de la m\'ethode des orbites de Kirillov aux groupes r\'esolubles exponentiels est due \`a P. Bernat \cite {BCD}. Soit $G$ un groupe de Lie r\'esoluble exponentiel connexe
et simplement connexe d'alg\`ebre de
Lie $\g g$, soit $f\in \g g^*$ et soit $\g h$ une polarisation
r\'eelle en $f$. Notons $H$ le sous-groupe $\exp \g h$. Soit
$O_f$ l'orbite coadjointe de  $f$ et $d$=dim $O_f$. Soit $\chi_f$ le caract\`ere de $H$ d\'efini par~:
$$\chi_f(\exp X)=e^{-i<f,X>}.\eqno{(\IV.1.1)}$$
Soit $\mop{pr}$ la surjection de $O_f$ sur $G/H$ d\'efinie par~:
$$\mop{pr}(g.f)=gH.\eqno{(\IV.1.2)}$$
Celle-ci est bien d\'efinie car $H$ contient le stabilisateur $G_f$ de $f$. On consid\`ere la restriction du crochet de Poisson de $\g g^*$ \`a la feuille symplectique $O_f$. Suivant \cite {Pe2} on d\'esigne par $\Cal E^0(O_f)$ l'espace des fonctions $\varphi=\psi\circ\mop{pr}$, o\`u $\psi\in C^\infty(G/H)$, et on d\'esigne par $\Cal E^1(O_f)$ le normalisateur de $\Cal E^0(O_f)$ dans $C^\infty(O_f)$ pour le crochet de Poisson. Tout $X\in\g g$ peut se voir comme une fonction lin\'eaire sur $\g g^*$. Elle induit par restriction une fonction $C^\infty$ sur $O_f$ que l'on note encore $X$.
\lemme{IV.1.1}
Pour tout $X\in\g g$ la fonction $X$ appartient \`a $\Cal E^1(O_f)$.
\dem
Le champ hamiltonien $H_X$ associ\'e \`a la fonction $X$ est \'egal au champ fondamental donn\'e par l'action coadjointe~:
$$H_X\varphi(\eta)={d\over dt}\restr{t=0}\varphi(\exp(-tX).\eta),
\ \eta\in O_f.
\eqno{(\IV.1.3)}$$
on a $\varphi\in \Cal E^0(O_f)$ si et seulement si $\varphi(gh.f)=\varphi(g.f)$ pour tout $g\in G$ et $h\in H$, ou encore si $H_X\varphi=0$ pour tout $X\in\g h$. Il est donc clair que pour tout $\varphi\in \Cal E^0(O_f)$ et pour tout $t\in \R$, $\eta\mapsto\varphi(\exp -tX.\eta)$ appartient aussi \`a $\Cal E^0(O_f)$. On en d\'eduit que $H_X.\varphi=\{X,\varphi\}$ appartient aussi \`a $\Cal E^0(O_f)$.
\qed
Soit $\tau=(\uple\tau{d/2}):G/H\fleche 6^\sim \Omega\subset\R^{d/2}$ une carte globale. On supposera ici que l'ouvert $\Omega$ est \'egal \`a $\R^{d/2}$ en entier. C'est toujours possible dans le cas d'un groupe r\'esoluble exponentiel (\cite {BCD}, \cite{AC}) en choisissant une base $(\uple X{d/2})$ coexponentielle \`a $\g h$ dans $\g g$ et en posant~:
$$\tau\inver(\uple x{d/2})=\exp{x_1X_1}\cdots\exp{x_{d/2}X_{d/2}}H.$$
On supposera par la suite que la carte globale $\tau$ est d\'efinie comme ci-dessus. Soient $(\uple q{d/2})$ les fonctions dans $\Cal E^0(O_f)$ d\'efinies par~:
$$q_j=\tau_j\circ\mop{pr}.\eqno{(\IV.1.4)}$$ Nous utiliserons les deux  r\'esultats suivants (\cite {Pe2} Th\'eor\`emes 2.2.2 et 3.2.3)~:
\th{IV.1.2 \rm (N. V. Pedersen)}
Il existe une famille $(\uple p{d/2})$ dans $\Cal E^1(O_f)$ telle que $(\uple p{d/2},\uple q{d/2})$ forme une carte de Darboux globale~:
$$\Phi~: O_f\fleche 6^\sim W_f\times\R^{d/2}$$
o\`u $W_f$ est un ouvert de $\R^{d/2}$, c'est-\`a-dire~:
$$\{p_i,p_j\}=\{q_i,q_j\}=0,\hskip 12mm \{p_i,q_j\}=\delta_i^j.\eqno{(\IV.1.5)}$$
 Cet ouvert est $\R^{d/2}$ en entier si et seulement si la polarisation $\g h$ v\'erifie la condition de Pukanszky~:
$$H.f=f+\g h^\perp.\eqno{(\IV.1.6)}$$
\ndem
\th{IV.1.3 \rm (N.V. Pedersen)}
Soit $(p,q)$ un syst\`eme de coordonn\'ees de Darboux globales~: $O_f\fleche 4^\sim W_f\times\R^{d/2}\subset \R^d$ comme dans le th\'eor\`eme pr\'ec\'edent. Alors~:
\smallskip
1°) dans ces coordonn\'ees l'espace $\Cal E^0(O_f)$ s'identifie \`a l'espace des fonctions qui ne d\'ependent que de $q$, et l'espace $\Cal E^1(O_f)$ s'identifie \`a l'espace des fonctions~:
$$\varphi(p,q)=\sum_{u=1}^{d/2}a_u(q)p_u+a_0(q)\eqno{(\IV.1.7)}$$
o\`u $a_0,\uple a{d/2}\in C^\infty(\R^{d/2})$. En particulier pour tout $X\in\g g$ la fonction associ\'ee lue dans la carte s'\'ecrit~:
$$X(p,q)=\sum_{u=1}^{d/2}a_{X,u}(q)p_u+a_{X,0}(q).\eqno{(\IV.1.8)}$$
o\`u les $a_{X,u},u=0,\ldots , d/2$ sont des fonctions dans $C^\infty(\R^{d/2})$.\smallskip
2°) Il existe une unique repr\'esentation unitaire fortement continue $\rho$ de $G$ dans $L^2(\R^{d/2})$ telle que l'espace $\Cal H_\rho^\infty$ de ses vecteurs $C^\infty$ contienne $C^\infty_c(\R^{d/2})$, et telle que pour tout $\xi\in~C^\infty(\R^{d/2})$ on a~:
$$\rho(X)\xi(t)=\sum_{u=1}^{d/2}a_{X,u}(t){\partial\xi(t)\over\partial t_u}-ia_{X,0}(t)\xi(t)+{1\over 2}\Bigl(\sum_{u=1}^{d/2}{\partial a_{X,u}\over\partial t_u}(t)\Bigr)\xi(t).\eqno{(\IV.1.9)}$$
Cette repr\'esentation est \'equivalente \`a l'induite $\mop{Ind}_H^G\chi_f.$
\ndem
L'expression des fonctions $X(p,q)$ d\'ecoule imm\'ediatement du lemme IV.1.1. Les fonctions $a_{X,u}$ sont de plus enti\`eres \`a croissance au plus exponentielle \cite{AC Th\'eor\`eme 1.6}. On peut tirer une information suppl\'ementaire sur la valeur en $q=0$ des fonctions $a_{X,u}$~:
\lemme{IV.1.4}
1) Pour tout $X\in\g g$ on a $a_{X,0}(0)=<f,X>$.
\smallskip
2) Pour tout $X\in\g h$ et pour tout $u=1,\ldots, d/2$ on a $a_{X,u}(0)=0$.

3) Pour tout $j=1,\ldots, d/2$ et pour tout $u=1,\ldots, d/2$ on a~:
$$a_{X_j,u}(0)=-\delta_j^u.$$
\dem
La premi\`ere assertion d\'ecoule imm\'ediatement du fait que le point de $O_f$ de coordonn\'ees $(0,0)$ est $f$. La deuxi\`eme assertion vient du fait que l'ensemble des points de coordonn\'ees $(p,q)$ avec $q=0$ est $H.f$, qui est inclus dans $f+\g h^\perp$. Quant au troisi\`eme point, il d\'ecoule d'un calcul direct~:
$$\eqalign{a_{X_j,u}(0)	&={\partial\over\partial p_u}X_j(p,0)	\cr
			&=-\{q_u,X_j\}(p,0)			\cr
			&=\{X_j,q_u\}(p,0)			\cr
			&={d\over dt}\restr{t=0}q_u\Bigl(
				\exp -tX_j.(p,0)\Bigr)		\cr
	&={d\over dt}\restr{t=0}q_u\Bigl(p,(0,\ldots,-t,\ldots 0)\Bigr)
		\hbox{ \qquad ($-t$ en position $j$)}\cr
			&=-\delta_j^u.				\cr}$$
\qed

\alinea{IV.2. Construction d'une repr\'esentation $\pi_\nu$ de $\Cal A$}
\qquad On applique la construction pr\'ec\'edente \`a une famille $(G_\hbar)_{\hbar\in\R-\{0\}}$ de groupes exponentiels. Plus pr\'ecis\'ement, soit $(\g g_\hbar)_{\hbar\in\R}$ la famille d'alg\`ebres de Lie r\'esolubles exponentielles d\'efinie par un m\^eme espace vectoriel sous-jacent $\g g$, avec le crochet~:
$$[X,Y]_\hbar=\hbar[X,Y].$$
On notera $\exp_\hbar$ l'exponentielle de $\g g_\hbar$ dans $G_\hbar$. Soit $f\in\g g^*$. L'orbite coadjointe $O_{f,\hbar}=G_\hbar.f\subset\g g^*$ est la m\^eme pour tout $\hbar\not =0$, mais la structure de Poisson d\'epend de $\hbar$~: pour tout $\eta\in\g g$ et $\varphi,\psi\in C^\infty(\g g^*)$ on a~:
$$\{\varphi,\psi\}_\hbar(\eta)=\hbar\{\varphi,\psi\}_1(\eta)
=\hbar<\eta,\,[d\varphi(\eta),d\psi(\eta)]>.\eqno{(\IV.2.1)}$$
On notera toujours $O_f$ cette orbite commune, lorsque la structure de Poisson n'a pas besoin d'\^etre pr\'ecis\'ee. En revanche l'orbite $O_{f,0}$ d\'eg\'en\`ere et se r\'eduit au point $f$, puisque le groupe $G_0$ est ab\'elien. 
\ssq
Il existe un sous-espace $\g h$ de $\g g$ qui est une polarisation r\'eelle de $f$ dans $\g g_\hbar$ pour tout $\hbar\not =0$. Soit $H_\hbar=\exp_\hbar \g h\subset G_\hbar$. Soit $\chi_{f,\hbar}$ le caract\`ere de $H_\hbar$ d\'efini par~:
$$\chi_{f,\hbar}(\exp_\hbar X)=e^{-i<f,X>}.\eqno{(\IV.2.2)}$$
On utilise les r\'esultats rappel\'es au \S\ IV.1 pour construire une r\'ealisation simultan\'ee de toutes les induites $\mop{Ind}_{H_\hbar}^{G_\hbar}\chi_{f,\hbar}$ dans le m\^eme espace $L^2(\R^{d/2})$. En effet si $(p,q)$ est le syst\`eme de coordonn\'ees de Darboux globales du th\'eor\`eme IV.1 pour $O_{f,1}$ (correspondant \`a $\hbar=1$), alors pour tout $\hbar\not =0$, un syst\`eme de coordonn\'ees de Darboux globales pour $O_{f,\hbar}$ est donn\'e par $(p,q')=(p,\hbar\inver q)=(\uple p{d/2},\uple{\hbar\inver q}{d/2})$. Pour tout $X\in\g g$ la fonction correspondante lue dans cette carte s'\'ecrit~:
$$\eqalign{X(p,q')	&=\sum_{u=1}^{d/2} a_{X,u}(q)p_u+a_{X,0}(q)\cr
&=\sum_{u=1}^{d/2} a_{X,u}(\hbar q')p_u+a_{X,0}(\hbar q').\cr}$$
Il existe donc pour tout $\hbar\not =0$, d'apr\`es le th\'eor\`eme IV.1.3, une unique repr\'esentation unitaire fortement continue $\rho_\hbar$ de $G_\hbar$ dans $L^2(\R^{d/2})$ telle que l'espace $\Cal H_{\rho_\hbar}^\infty$ de ses vecteurs $C^\infty$ contienne $C^\infty_c(\R^{d/2})$, et telle que~:
$$\rho_\hbar(X)\xi(t)=\sum_{u=1}^{d/2}a_{X,u}(\hbar t){\partial\xi(t)\over\partial t_u}-ia_{X,0}(\hbar t)\xi(t)+{\hbar\over 2}\Bigl(\sum_{u=1}^{d/2}{\partial a_{X,u}\over\partial t_u}(\hbar t)\Bigr)\xi(t).\eqno{(\IV.2.3)}$$
Cette repr\'esentation est \'equivalente \`a l'induite $\mop{Ind}_{H_\hbar}^{G_\hbar}\chi_{f,\hbar}$. Elle est irr\'eductible si et seulement si la polarisation $\g h$ v\'erifie la condition de Pukanszky $H_\hbar.f=f+\g h^\perp$ pour un quelconque $\hbar\not =0$.
\ssq
On peut par ailleurs facilement r\'ealiser l'induite $\mop{Ind}_{H_0}^{G_0}\chi_{f,0}$ dans $L^2(\R^{d/2})$, ce qui donne~:
$$\rho_0(\exp tX_j).\varphi(\uple t{d/2})
	=\varphi(t_1,\ldots ,t_{j-1},t_j-t,t_{j+1},\ldots ,t_{d/2}),$$
et pour tout $X\in\g h$~:
$$\rho_0(\exp tX).\varphi(\uple t{d/2})
	=e^{-it<f,X>}\varphi(\uple t{d/2}).$$
En diff\'erenciant en $t=0$ on obtient~:
$$\rho_0(X_j).\varphi(\uple t{d/2})=-{\partial\over\partial t_j}\varphi(\uple t{d/2}),\eqno{(\IV.2.4)}$$
et pour tout $X\in\g h$~:
$$\rho_0(X).\varphi(\uple t{d/2})=-i<f,X>.\varphi(\uple t{d/2}).\eqno{(\IV.2.5)}$$
Gr\^ace au lemme IV.1.4 les notations sont coh\'erentes~: la repr\'esentation $\rho_0$ de $G_0$ dans $C^\infty_c(\R^{d/2})$ s'obtient en prolongeant \`a $\hbar=0$ la formule (IV.2.3).
\ssq
Soit maintenant $K$ un compact fix\'e de $\R^{d/2}$ d'int\'erieur non vide, et soit $M=C^\infty_K(\R^{d/2})$ l'espace des fonctions $C^\infty$ \`a support dans $K$. Cet espace est pour tout r\'eel $\hbar$ un sous-module du module des vecteurs $C^\infty$ de $\rho_\hbar$. On munit $M$ de la topologie de Fr\'echet d\'efinie par les seminormes~:
$$N_k(\varphi)=\mopl{sup}_{|\alpha|\le k}\mopl{sup}_{t\in K}|D^\alpha\varphi(t)|.$$
On reprend les notations du \S\ II. Le dual topologique de $M$ est l'espace des distributions de support inclus dans $K$. Soient $\varphi\in M$ et $T\in M'$. Comme les fonctions $a_{X,u}$ sont enti\`eres, on voit en appliquant le lemme III.0.1 que l'expression $<T,\,\rho_\hbar(w)\varphi>$ est enti\`ere en $\hbar$ pour tout $w\in S(\g g)$, o\`u l'on a identifi\'e $S(\g g)$ et $\Cal U(\g g_\hbar)$ via l'isomorphisme de Duflo. En posant~:
$$\pi_{\nu_0}(X)=i\rho_{-i\nu_0}(X)\eqno{(\IV.2.6)}$$
on obtient une famille $(\pi_{\nu_0})$ de repr\'esentations de $(A,*_{\nu_0})$ dans $M$ v\'erifiant les conditions de la proposition II.2.1 (avec $R=+\infty$). En appliquant la proposition II.2.1 et le th\'eor\`eme II.3.1 on a donc structure de module topologiquement libre faiblement convergent de rayon infini fortement unitaire sur $\Cal M=M[[\nu]]$. L'expression
de la repr\'esentation $\pi_\nu$ de l'alg\`ebre d\'eform\'ee $\Cal A$ dans $\Cal M$ s'obtient comme suit~: pour tout $X\in\g g$ on a~:
$$\pi_\nu (X)\xi (t) =\sum_{u=1}^{d\over 2}i a_{X,u}(-i\nu t)
{{\p \xi}\over {\p t_u}}(t)+a_{X,0}(-i\nu t)\xi (t)
 +{\nu\over 2} \sum_{u=1}^{d\over 2}
{{\p \over {\p t_u}}} a_{X,u}(-i\nu t) \xi (t),\eqno{(\IV.2.7)}$$
et l'action d'un \'el\'ement de $\Cal A$ s'obtient via l'identification $\Cal A=\Cal U_\nu(\g g^\C)$ donn\'ee par l'isomor\-phisme de Duflo.
\alinea{IV.3. Calcul de la vari\'et\'e caract\'eristique}
\th{IV.3.1}
Soient $\g g$ une alg\`ebre de Lie r\'esoluble exponentielle, $f$ un \'el\'ement de $\g g^*$
et $\g h$ une polarisation r\'eelle en $f$. Soit $d$ la dimension de l'orbite coadjointe de $f$, soit $K$ un compact de $\R^{d/2}$ d'int\'erieur non vide. Soit $\pi_\nu$ la repr\'esentation de l'alg\`ebre d\'eform\'ee ${\Cal A}$ dans le module topologiquement libre convergent $M=C^\infty_K(\R^{d/2})[[\nu]]$ construite au paragraphe IV.2 \`a partir de ces donn\'ees. Alors,
$$V(\pi_\nu)=f+\g h^\bot.$$
En particulier, si $\g h$ v\'erifie la condition de Pukanszky on a $V(\pi_\nu)=H.f$.
\dem
Faisant $\nu=0$ dans (IV.2.7), on obtient~:
$$\pi_0 (X) =\sum_{u=1}^{d\over 2}i a_{X,u}(0){{\p }\over {\p t_u}}+a_{X,0}(0).$$
D'apr\`es le lemme IV.1.4 on a donc~:
$$\cases{\pi_0(X_j)=-i{\partial\over\partial t_j}, j=1,\ldots ,d/2	\cr
	\pi_0(X)=<f,X> \hbox{ pour tout } X\in\g h.\cr}$$

L'annulateur de $\pi_0$ est l'id\'eal de $S(\g g)$
 engendr\'e par les $X-<f,X>,\,X\in\g h$, d'o\`u l'\'egalit\'e~:
$$V(\pi_\nu)=f+\g h^\bot.$$
\qed
\alinea{IV.4. Le cas nilpotent}
Dans les exemples III.1 et III.2 la vari\'et\'e de Poisson
caract\'eristique co\"\i ncide avec l'orbite coadjointe. Cette
propri\'et\'e se g\'en\'eralise \`a toute repr\'esentation
unitaire irr\'eductible d'un groupe nilpotent~: soit $\g g$ une
alg\`ebre de Lie r\'eelle nilpotente de dimension $n$. Soit $f\in\g g^*$, soit $\g h$ une polarisation r\'eelle en $f$ et soit $\rho_\hbar$ la repr\'esentation unitaire du
groupe simplement connexe $G_\hbar =\exp \g g_\hbar$ donn\'ee par la construction du \S\ IV.2. Elle est irr\'eductible car $\g h$ v\'erifie toujours la condition de Pukanszky dans le cas nilpotent. L'orbite coadjointe
$O_f\subset\g g^*$ est une sous-vari\'et\'e alg\'ebrique, donn\'ee
par l'annulation de polyn\^omes \`a valeurs r\'eelles
$(Q_j)_{j=1,\ldots , n-d}$ o\`u $d$ d\'esigne la dimension de l'orbite.
\ssq
 Soit $\tau$ l'isomorphisme
de Duflo, qui se r\'eduit ici \`a la sym\'etrisation. D'apr\`es le
th\'eor\`eme 2.3.2 de \cite{Pe1} (adapt\'e ici \`a nos conventions
de signe dans la d\'efinition du caract\`ere $\chi_f$),
l'annulateur de $\rho_\hbar$ dans $\Cal U(\g g_\hbar)$ est
l'id\'eal engendr\'e par les $u_j$ avec
$\tau\inver(u_j)(\eta)=Q_j(i\eta)$. On voit alors que l'annulateur de la repr\'esentation $\pi_\nu$
 de l'alg\`ebre
 d\'eform\'ee $\Cal A$ est
engendr\'e par les $Q_j$. L'id\'eal
$\mop{Ann}\pi_\Omega^\nu/(\mop{Ann}\pi_\Omega^\nu\cap\nu\Cal A)$ de $S(\g g)$
est donc \'egalement engendr\'e par les $Q _j$, d'o\`u
l'\'egalit\'e entre l'orbite $\Omega$ et la vari\'et\'e de Poisson
caract\'eristique. On a donc~:
\th{IV.4.1} Avec les notations
ci-dessus,
$$V\!A(\pi_\nu)=\Omega.$$
\ndem
{\sl Remarque\/}~: ce th\'eor\`eme est en g\'en\'eral faux pour un groupe r\'esoluble exponentiel non nilpotent, comme le montrent les exemples III.3 et III.4. Dans ces deux exemples la vari\'et\'e de Poisson caract\'eristique est \'egale \`a la fermeture de Zariski de l'orbite.
\paragraphe{R\'ef\'erences~:}
\vskip -3mm
\bib{A}S. Araki, {On root systems and an infinitesimal classification of irreducible symmetric
spaces\/}, J. of Math. Osaka City Univ. 13, 1-34 (1962).
\bib{AC}D. Arnal, J-C. Cortet, {\sl Repr\'esentations $*$ des groupes exponentiels\/}, J. Funct. Anal. 92-1, 103-135 (1990).
\bib{AM}M. Andler, D. Manchon, {\sl Op\'erateurs aux diff\'erences
finies, calcul pseudo-diff\'erentiel et repr\'esentations des groupes de
Lie}, J. Geom. Phys. 27, 1-29 (1998).
\bib{ABLS}D. Arnal, A. Baklouti, J. Ludwig, M. Selmi, {\sl
Separation of unitary representations of exponential Lie
groups\/}, J. Lie Theory, vol. 10, 399-410 (2000).
\bib{ABM}D. Arnal, N. Ben Amar, M. Masmoudi, {\sl Cohomology of good
graphs and Kontsevich linear star products}, Lett. Math. Phys. 48, 291-
306 (1999).
\bib{AMM}D. Arnal, D. Manchon, M. Masmoudi, {\sl Choix des signes
pour la formalit\'e de M. Kontsevich}, math.QA/0003003, Pac. J. Math.
\bib{BBR}A. Baklouti, C. Benson and G. Ratcliff,
{\sl Moment sets and the unitary dual of a nilpotent Lie group},
J. Lie Theory, Vol 11, 135-154 (2001).
\bib{BCD}P. Bernat, N. Conze, M. Duflo et al. {\sl Repr\'esentations des
groupes de Lie r\'esolubles}, Monographies de la Soc. Math. France No 4,
Dunod, Paris 1972.
\bib{BGHHW}M. Bordemann, G. Ginot, G. Halbout, H-C. Herbig, S. Waldmann,
{\sl Star-repr\'esentations sur des sous-vari\'et\'es co-isotropes\/}, math.QA/0309321 (2003).
\bib{BW}M. Bordemann, S. Waldmann, {\sl Formal GNS construction and
states in deformation quantization}, Comm. Math. Phys. 195, 549-583
(1998).
\bib{Bu-W}H. Bursztyn, S. Waldmann, {\sl $*$-ideals and formal Morita
equivalence of $*$-algebras}, Int. J. Math. 12 No 5, 555-577 (2001).
\bib {B}N. Bourbaki, {\sl Alg\`ebre commutative}, Hermann, Paris 1961.
\bib{CF1}A. Cattaneo, G. Felder, {\sl A path integral approach to the Kontsevich quantization
formula}, Comm. Math. Phys. 212 N°3, 591-611 (2000).
\bib{CF2}A. Cattaneo, G. Felder, {\sl Coisotropic submanifolds in Poisson geometry and branes
in the Poisson sigma model\/}, math.QA/0309180 (2003).
\bib{CFT}A. Cattaneo, G. Felder, L. Tomassini, {\sl From local to global
deformation quantization of Poisson manifolds}, math.QA/0012228 (2000).
\bib{Di}G. Dito, {\sl Kontsevich star product on the dual of a Lie
algebra}, Lett. Math. Phys. 48, 307-322 (1999).
\bib{Dix}J. Dixmier, {\sl Alg\`ebres enveloppantes}, Gautier-Villars, Paris 1974.
\bib{D1}M. Duflo, {\sl Op\'erateurs diff\'erentiels bi-invariants sur un
groupe de Lie, \/} Ann. Sci. Ec. Norm. Sup. 4e s\'erie t. 10, 265-288 (1977).
\bib{D2}M. Duflo, {\sl Sur la classification des id\'eaux primitifs dans
l'alg\`ebre enveloppante d'une alg\`ebre de Lie semi-simple},
Ann. of Math. (2) 105, no. 1, 107-120 (1977).
\bib{E}B. Enriquez, {\sl Quantization of Lie bialgebras and shuffle algebras of Lie algebras\/}, Selecta Math. 7, no. 3, 321-407 (2001).
\bib{EHW}T. Enright, R. Howe, N. Wallach, {\sl A classification of unitary highest
weight modules\/}, Progr. Math. 40, 97-143,(1983).
\bib{EK}P. Etingof, D. Kazhdan, {\sl Quantization of Lie bialgebras, I},
Selecta Math. 2, 1-41 (1996).
\bib{FS}G. Felder, B. Shoikhet, {\sl Deformation quantization with traces},
math.QA/0002057.
\bib {G}O. Gabber, {\sl The integrability of the characteristic variety},
Am. J. Math. 103, No 3, 445-468 (1981).
\bib {GM}M. Granger, Ph. Maisonobe, {\sl A basic course on differential
modules}, in {\sl $\Cal D$-modules coh\'e\-rents et holonomes}, les
cours du CIMPA, Hermann, Paris 1993.
\bib {Gui}V. Guinzburg, {\sl On primitive ideals}, math.RT/0202079.
\bib{H}S. Helgason, {\sl Differential geometry and symmetric spaces}, Academic Press, 1962.
\bib{J}A. Joseph, {\sl On the classification of primitive ideals in the
enveloping algebra of a semisimple Lie algebra}, Lect. Notes in Math.
1024 , 30-78.
\bib{Ja}J.C. Jantzen, {\sl Einh\"ullende Algebren halbeinfacher
Lie-Algebren}, Springer (1983).
\bib{Jak}H.P. Jakobsen, {\sl Hermitian symmetric spaces and their unitary highest weight
 modules}, J. Funct. Anal. 52 No 3, 385-412 (1983).
\bib{Ka}Ch. Kassel, {\sl Quantum Groups\/}, Springer (1995).
\bib{Ki}A. A. Kirillov, {\sl Elements of the Theory of Representations},
Springer (1976).
\bib {K}M. Kontsevich, {\sl Deformation quantization of Poisson
manifolds I,\/} math.QA/9709040.
\bib {MT}D. Manchon, Ch. Torossian, {\sl Cohomologie tangente et
cup-produit pour la quantification de Kontsevich}, math.QA/0106205.
\bib {Pe1}N. V. Pedersen, {\sl On the infinitesimal kernel of irreducible
representations of nilpotent Lie groups}, Bull. Soc. Math. France 112, 42
3-467 (1984).
\bib {Pe2}N. V. Pedersen, {\sl On the symplectic structure of coadjoint orbits
of (solvable) Lie groups and applications. I}, Math. Ann. 281, 633-669 (1988).
\bib{Sha}N.N. Shapovalov, {\sl A certain bilinear form on the universal enveloping algebra
of a complex semisimple Lie algebra\/} Funct. Anal. Appl. 6 , no. 4, 307-312 (1972).
\bib {Sh}B. Shoikhet, {\sl Vanishing of the Kontsevich integrals of the
wheels}, math.QA/0007080.
\bib{T} D. Tamarkin, {Another proof of M. Kontsevich formality theorem\/}, math.QA/9803025 (1998).
\bib{Tr}F. Treves, {\sl Topological vector spaces, distributions and kernels\/},
Academic Press, New York (1967).
\bib {V}P. Vanhaecke, {\sl Integrable systems in the realm of algebraic
geometry}, Lect. Notes in Math. No 1638, Springer (1996).
\bib{Wa}G. Warner, {\sl Harmonic analysis on semi-simple Lie groups I},
Grundlehren der mathematischen Wissenschaften No 188, Springer (1972).
\bib {W}A. Weinstein, {\sl The local structure of Poisson manifolds},
J. Diff. Geom. 18, 523-557 (1983).
\bib{Y}A. Yekutieli, {\sl On deformation quantization in algebraic geometry\/}, math.AG/0310399 (2003).
\bye